\documentclass[aop,MSNbibl,citesort,dvips]{arximspdf}
\usepackage{mathbh}
\usepackage{newsym2e}
\usepackage{pst-node}
\usepackage{graphicx}

\doi{10.1214/10-AOP561}
\volume{39}
\issue{4}
\pubyear{2011}
\firstpage{1205}
\lastpage{1242}

\makeatletter

\makeatletter
\let\color@gray\@gobbletwo
\makeatother

\newcommand{\eqref}[1]{(\ref{#1})}

\newcommand{\implies}{\Longrightarrow}

\newcommand{\xrightarrow}{\mathop{\longrightarrow}\limits_}
\newcommand{\xxrightarrow}{\limits^}

\newtheorem{theorem}{Theorem}[section]
\newtheorem{lemma}[theorem]{Lemma}
\newtheorem{coro}[theorem]{Corollary}
\newtheorem{prop}[theorem]{Proposition}
\newtheorem{conj}[theorem]{Conjecture}

\newproclaim{example}[theorem]{Example}
\newproclaim{fact}[theorem]{Fact}
\newproclaim{defn}[theorem]{Definition}

\newcommand{\eps}{\varepsilon}
\renewcommand{\P}{\mathbb{P}}
\newcommand{\E}{\mathbb{E}}
\newcommand{\R}{\mathbb{R}}
\newcommand{\Z}{\mathbb{Z}}

\newcommand{\cM}{\mathcal{M}}

\newcommand{\F}{\mathcal{F}}
\newcommand{\eqd}{\stackrel{d}{=}}

\newcommand{\bx}{\mathbf{x}}
\newcommand{\bB}{\mathbf{B}}
\newcommand{\muI}{\mu^{(2)}}
\newcommand{\muII}{\tilde\mu^{(2)}}
\newcommand{\ind}{\mathbh{1}}

\renewcommand{\bar}[1]{{\overline{#1}}}

\newcommand{\bin}{\operatorname{Bin}}
\newcommand{\geom}{\operatorname{Geom}}

\newcommand{\id}{\operatorname{id}}
\newcommand{\Leb}{\operatorname{{\mathfrak Leb}}}

\newenvironment{configarray}[1] {%
\begin{array}[c]{*{#1}{ @{}c } @{}}%
}{%
\end{array}%
}

\newcommand{\config}[2]{\fbox{$\begin{configarray}{#1} #2 \end{configarray}$}}

\newcommand{\C}{\mbox{\gray{\CIW}}}
\newcommand{\B}{\raisebox{0.9pt}{\mbox{{$\circledast$}}}}
\newcommand{\D}{\mbox{\CIB}}

\makeatother

\begin{document}
\begin{frontmatter}

\title{The TASEP speed process\thanksref{T1}}
\runtitle{The TASEP speed process}

\thankstext{T1}{Research was carried out while all authors were at the
University of Toronto. Supported in part by NSERC.}

\begin{aug}
\author[A]{\fnms{Gideon} \snm{Amir}\corref{}\ead[label=e1]{amirgi@math.biu.ac.il}\ead[label=u1,url]{http://u.math.biu.ac.il/\textasciitilde amirgi/}},
\author[B]{\fnms{Omer} \snm{Angel}\ead[label=e2]{angel@math.ubc.ca}\ead[label=u2,url]{http://www.math.ubc.ca/\textasciitilde angel}} and
\author[C]{\fnms{Benedek} \snm{Valk\'{o}}\ead[label=e3]{valko@math.wisc.edu}\ead[label=u3,url]{http://www.math.wisc.edu/\textasciitilde valko}}

\runauthor{G. Amir, O. Angel and B. Valk\'{o}}
\affiliation{Bar Ilan University, University of British Columbia and
University~of~Wisconsin}
\address[A]{G. Amir \\
Department of Mathematics\\
Bar Ilan University\\
52900 Ramat Gan\\
Israel\\
\printead{e1}\\
\printead{u1}}
\address[B]{O. Angel \\
Department of Mathematics\\
University of British Columbia\\
Vancouver BC, V6T 1Z2\\
Canada\\
\printead{e2}\\
\printead{u2}}
\address[C]{B. Valk\'{o} \\
Department of Mathematics\\
University of Wisconsin--Madison\\
480 Lincoln Dr.\\
Madison, Wisconsin 53706\\
USA\\
\printead{e3}\\
\printead{u3}}
\end{aug}

\received{\smonth{8} \syear{2009}}
\revised{\smonth{4} \syear{2010}}

%
\begin{abstract}
In the multi-type totally asymmetric simple exclusion process (TASEP) on
the line, each site of $\Z$ is occupied by a particle labeled with some
number, and two neighboring particles are interchanged at rate one if
their labels are in increasing order. Consider the process with the
initial configuration where each particle is labeled by its position. It
is known that in this case a.s. each particle has an asymptotic speed
which is distributed uniformly on $[-1,1]$. We study the joint
distribution of these speeds: the \textit{TASEP speed process}.

We prove that the TASEP speed process is stationary with respect to the
multi-type TASEP dynamics. Consequently, every ergodic stationary measure
is given as a projection of the speed process measure. This generalizes
previous descriptions restricted to finitely many classes.

By combining this result with known stationary measures for\break TASEPs with
finitely many types, we compute several marginals of the speed process,
including the joint density of two and three consecutive speeds. One
striking property of the distribution is that two speeds are equal with
positive probability and for any given particle there are infinitely many
others with the same speed.

We also study the partially asymmetric simple exclusion process (ASEP).
We prove that the states of the ASEP with the above initial
configuration, seen as permutations of $\Z$, are symmetric in
distribution. This allows us to extend some of our results, including the
stationarity and description of all ergodic stationary measures, also to
the ASEP.
\end{abstract}

%

\end{frontmatter}

\section{Introduction} \label{S:intro}

The exclusion process on a graph describes a system of particles performing
continuous time random walks, interacting with other particles via
exclusion: attempted\vadjust{\goodbreak} jumps to occupied sites are suppressed. When the graph
is $\Z$ and particles jump only to the right at rate one the process is
called the \textit{totally asymmetric simple exclusion process}
(TASEP). We
denote configurations with $\eta\in\{1,\infty\}^{\Z}$ where particles are
denoted by $1$ and empty sites by $\infty$.\setcounter
{footnote}{1}\footnote{The common
practice is
to denote empty sites by $0$. However, under various common extensions of
the TASEP including those used here, it is convenient to denote empty
sites by a label larger than the labels of all particles.} The TASEP
is a
Markov process with generator
%
%
\begin{equation}
\label{eq:gen1}
L f(\eta)=\sum_{n} f(\sigma_n\eta) - f(\eta),
\end{equation}
where $\sigma_n$ is the operation that sorts the coordinates at $n,
n+1$ in
decreasing order
%
%
\begin{eqnarray}
(\sigma_n \eta)_n &=&\max(\eta_n,\eta_{n+1}),\qquad
(\sigma_{n} \eta)_{n+1}=\min(\eta_n,\eta_{n+1}),\nonumber\\[-8pt]\\[-8pt]
(\sigma_n \eta)_k&=&\eta_k \qquad\mbox{if $k\notin{n,n+1}$}.\nonumber
\end{eqnarray}

A \textit{second class particle} is an extra particle in the system
trying to
perform the same random walk while being treated by the normal (first
class) particles as an empty site. It is an intermediate state between a
particle and an empty site, and is denoted by a $2$.\footnote{$k$th class
particles will be denoted by $k$, even for $k=0$. That is why it is
convenient to use $\infty$ for holes rather than $0$.} This means that
the second class particle will jump to the left if there is a first class
particle there who decides to jump onto the second class particle. This is
still a Markov process, with the same generator \eqref{eq:gen1} and state
space $\{1,2,\infty\}^\Z$. Note that empty sites can just be considered as
particles with the highest possible class. Thus we can equally well
consider state space $\{1,2,3\}^\Z$ with holes represented by $3$'s.

More generally, we shall consider the multi-type TASEP which has the same
generator with state space $\R^\Z$. Thus we allow particle classes to be
nonintegers or negative numbers. If there are particles with maximal class
they can be considered to be holes. A special case is the $N$-type TASEP
(without holes) where all particles have classes in $\{1,\ldots,N\}$. If
particles of class $N$ are interpreted as holes instead of maximally
classed particles, this process becomes the traditional $(N-1)$-type TASEP
(with holes). To avoid confusion, from here on all multi-type
configurations shall be without holes. (Holes will appear only in
individual lines in the multi-line configurations defined below.)

The following result is this paper's foundation. We let $Y(t)$ denote the
TASEP configuration at time $t$, with $Y_n(t)$ the value at position $n$.
This strengthens results of Ferrari and Kipnis~\cite{FerrariKipnis} that
get the same limit in distribution.
\begin{theorem}[(Mountford and Guiol~\cite{MountfordGuiol})]\label{T:MG}
Consider the TASEP with initial condition
\[
Y_n(0) = \cases{ 1, &\quad$n<0$, \cr
2, &\quad$n=0$, \cr
3, &\quad$n>0$.}
\]
Let $X(t)$ denote the position of the second class particle at time $t$,
defined by $Y_{X(t)}(t) = 2$. Then $\frac{X(t)}{t}
\xrightarrow{t\to\infty}\xxrightarrow{a.s.} U$, where $U$ is a uniform
random variable on
$[-1,1]$.
\end{theorem}

Thus a second class particle with first class particles to its left and
third class particles to its right ``chooses'' a speed $U$, uniform in
$[-1,1]$ and follows that speed: $X(t) \sim U t$. (See~\cite{FP,FMP} for
alternative proofs of Theorem~\ref{T:MG}.)

Now, consider any other starting configuration such that $Y_n(0)<Y_0(0)$
for all $n<0$ and $Y_n(0)>Y_0(0)$ for all $n>0$. The particle starting
at 0
does not distinguish between higher classes, or between lower classes, so
its trajectory has the same law. This applies in particular to every
particle in a multi-type TASEP $Y$ with starting configuration $Y_n(0)=n$.
Let $X_n(t)$ be the location of particle $n$ at time $t$, so that
$Y_{X_n(t)}(t)=n$ [$X(t)$ is the inverse permutation of $Y(t)$]. An
immediate consequence is the following:
\begin{coro}[(The speed process)]\label{C:speed}
In the TASEP with starting configuration $Y_n(0)=n$, a.s. every particle
has a speed: for every $n$
\[
\frac{X_n(t)-n}{t} \xrightarrow{t\to\infty}\xxrightarrow{a.s.}
U_n,
\]
where $\{U_n\}_{n\in\Z}$ is a family of random variables, each uniform on
$[-1,1]$.
\end{coro}
\begin{defn}
The process $\{U_n\}_{n\in\Z}$ is called the \textit{TASEP speed process}.
Its distribution is denoted by $\mu$.
\end{defn}

\begin{figure}

\includegraphics{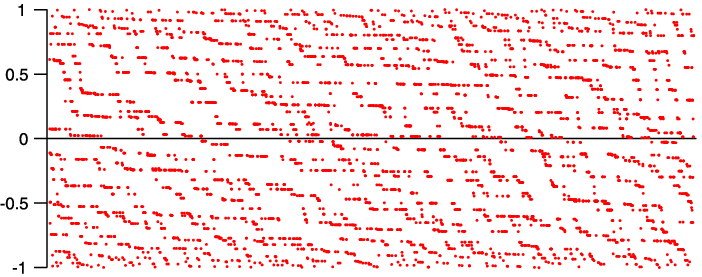}

\caption{The speed process: simulation of $U_n$ for $1\leq n \leq5000$,
from a simulation run to time $700\mbox{,}000$.}
\label{fig:speeds}
\end{figure}

Thus $\mu$ is a measure supported on $[-1,1]^\Z$. It is clear from
simulations (and our results below) that $\mu$ is not a product measure,
that is, that the speeds are not independent. Figure~\ref{fig:speeds} shows
%
%
a portion of the process. Some aspects of this process were studied in
\cite{FGM}.

\subsection{Main results}

In order to study the TASEP speed process we prove two results, which are
our main tools in understanding the joint distributions of speeds. These
results are of significant interest in and of themselves. The following is
a new and surprising symmetry of the TASEP. A version of this theorem was
proved in~\cite{DPRW}, in the context of the TASEP on finite intervals. We
extend it here also to the ASEP\footnote{Some sources use PASEP/ASEP,
respectively, for what other sources call ASEP/TASEP (PASEP stands for
partially$\ldots$). We adopt the latter convention.} (defined in
Section~\ref{SS:asep}).
\begin{theorem}\label{T:symmetry}
Consider the starting configuration $Y_n(0)=n$ and $X_n(t)$ as above. For
any fixed $t>0$ the process $\{X_n(t)\}_{n\in\Z}$ has the same
distribution as $\{Y_n(t)\}_{n\in\Z}$. This holds also for the ASEP.
\end{theorem}

At any time $t$ we have that $X(t)$ and $Y(t)$ are permutations of
$\Z$, one the inverse of the other. Thus this theorem implies that $Y(t)$
as a permutation has the same law as its inverse. It is not hard to see
that this holds only for a fixed $t$, and not as processes in $t$
[e.g.,
$X_0(t)$ changes by at most 1 at each jump].

The next result gives additional motivation for considering the speed
process, as it relates its law $\mu$ to stationary measures of the
multi-type TASEP (and ASEP).
\begin{theorem}\label{T:stationary}
$\mu$ is itself a stationary measure for the TASEP: the unique ergodic
stationary measure which has marginals uniform on $[-1,1]$.
\end{theorem}

This means that if we consider a TASEP in
$[-1,1]^\Z$ where
the initial configuration $Y(0)$ has distribution $\mu$ then at any
time $t$ the distribution of $Y(t)$
is also given by $\mu$.

It is known that the $N$-type process has ergodic stationary measures, and
that the distribution of $Y_n$ among the classes determines this
distribution uniquely. Standard techniques (see below) can be used to show
that the same holds also with infinitely many classes. Specifically, for
any distribution on $\R$ there is a unique ergodic stationary measure for
the TASEP with $Y_0$ (and any $Y_n$) having that distribution.
For any two nonatomic distributions on $\R$, these measures are
related by
applying pointwise a nondecreasing function to the particle classes (see
Lem\-ma~\ref{L:projection}), so every such measure can be deduced from
the measure
with marginals uniform on $[-1,1]$. If a distribution has atoms, then the
corresponding stationary measure can still be deduced from the speed
process' law $\mu$ in the same way, but the operation is nonreversible.
Thus we have the following characterization:
\begin{coro}\label{C:characterization}
Every ergodic stationary measure for the TASEP can be deduced from $\mu$
by taking the law of $\{F(U_n)\}_{n\in\Z}$ for some nondecreasing
function $F\dvtx[-1,1]\to\R$.
\end{coro}

%
\subsection{Results: Joint distribution}

\begin{figure}

\includegraphics{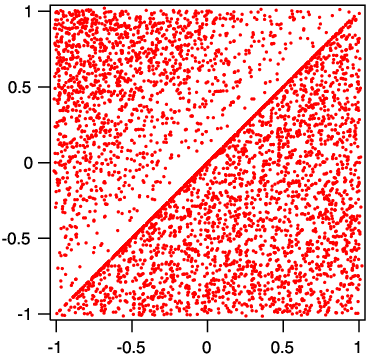}

\caption{The joint distribution of $U_0,U_1$: based on 5000 pairs from a
simulation run to time $25\mbox{,}000$.}
\label{fig:joint2}
\end{figure}

Computer simulations suggested early on that $U_0,U_1$ are not independent
(see Figure~\ref{fig:joint2}). Recent results of Ferrari, Goncalves and
Martin~\cite{FGM} confirm this prediction. They proved (among other things)
%
%
that the probability that particle $0$ eventually overtakes particle $1$
(we identify a particle with its class) is $2/3$. It follows that $\P(U_0
\geq U_1) \geq2/3$ (not necessarily equal since $U_0=U_1$ does not a
priori imply overtaking). Our first theorem describing the joint
distribution of speeds is the following:
\begin{theorem}\label{T:joint2}
The joint distribution of $(U_0,U_1)$, supported on $[-1,1]^2$, is
\[
f(x,y) \,dx \,dy + g(x) \ind(x=y) \,dx
\]
with
\[
f(x,y) = \cases{
\dfrac{1}{4}, &\quad$x>y$,\vspace*{2pt}\cr
\dfrac{y-x}{4}, &\quad$x \leq y$,}\qquad
g(x) = \frac{1-x^2}{8}.
\]
In particular, $\P(U_0>U_1) = 1/2$, $\P(U_0=U_1)=1/6$ and
$\P(U_0<U_1)=1/3$.
\end{theorem}

\subsubsection*{Remarks}
Note that the density in $\{U_0<U_1\}$ (linear in $U_1-U_0$, so that there
is repulsion between the speeds) can be deduced using only
Theorem~\ref{T:symmetry} (we do not include this argument here).
However, proving
the---seemingly simpler---constant density on $\{U_0>U_1\}$ and
deriving the singular component on the diagonal requires the power of
Theorem~\ref{T:stationary}. It is interesting to compare the power of
Theorem~\ref{T:symmetry} with that of the methods of~\cite{FGM}. It
appears that
both methods run into similar difficulties and have similar consequences,
suggesting a fundamental connection (there are also some parallels in the
proofs). Specifically, can the density in the region $\{U_0<U_1\}$ be
derived using the techniques of~\cite{FMP}? Finally, it is interesting that
our proof relies nontrivially on the extension of the TASEP to infinitely
many different classes of particles, though the question and answer can
both be posed using only $4$ classes (including holes). A similar remark
holds about some other results below as well.

Additional information about the joint distribution of speeds is
derived in
Section~\ref{S:mult}. We derive certain properties of the $n$-dimensional
marginals of $\mu$, and in Theorem~\ref{T:3speeds} we compute the joint
distribution of three consecutive speeds.

A surprising aspect of Theorem~\ref{T:joint2} is that there is a positive
probability ($1/6$) that $U_0=U_1$, even though each is uniform on
$[-1,1]$. Indeed, for any two particles there is a positive probability
that their speeds are equal. This phenomenon can be thought of as a
spontaneous formation of ``convoys,'' sets of particles that have the same
asymptotic speed, so their trajectories remain close. Our next result gives
a full description of such a convoy.
\begin{theorem}\label{T:convoys}
Let the convoy of $0$ be $C_0 = \{j \dvtx U_j=U_0\}$, that is, the set
of all
$j$ with the same speed as $0$. Then $C_0$ is $\mu$-a.s. infinite with
$0$ density. Moreover, conditioned on $U_0$, $C_0$ is a renewal process,
and the nonnegative elements of $C_0$ have the same law as the times of
last increase of a random walk conditioned to remain positive, with step
distribution
\[
\P(X=1) = \P(X=-1) = \frac{1-U_0^2}{4},\qquad
\P(X=0) = \frac{1+U_0^2}{2}.
\]
\end{theorem}

The ``times of last increase'' of a walk $Z$ are those indices $n$ for
which $m>n$ implies $Z_m>Z_n$. In particular the convoys are infinite and
they provide a translation invariant partition of the integers into
infinitely many infinite sets. The convoys are essentially the process with
$0$ density for second class particles, seen from a second class particle,
as studied by Ferrari, Fontes and Kohayakawa in~\cite{FFK}.

\subsection{The ASEP}\label{SS:asep}

As the name suggests, the totally asymmetric simple exclusion process
is an
extremal case of the \textit{asymmetric simple exclusion process}: the ASEP.
The ASEP is defined in\vadjust{\goodbreak} terms of a parameter $p\in(1/2,1]$, with $p=1$ being
the TASEP. While most quantities involved depend on $p$, the dependence
will usually be implicit.

In the ASEP particles jump one site to the right at rate $p\in(1/2,1]$ and
to the left at rate $\bar{p}=1-p$ (we use the convention $\bar{x}=1-x$).
The generator of this Markov process is
%
%
\begin{equation}\label{eq:gen_asep}
L f(\eta)= \sum_{n} p \bigl(f(\sigma_n\eta) - f(\eta)\bigr) + \bar{p}
\bigl(f(\sigma_n^*\eta) - f(\eta)\bigr),
\end{equation}
where $\sigma_n$ and $\sigma_n^*$ sort the values in $n, n+1$ in decreasing
and increasing order, respectively.

While some of the questions above make sense also in this setting,
there is
a key difficulty in that the analogue of Theorem~\ref{T:MG} for the
ASEP (conjectured below) is still unproved. Using the methods of
Ferrari and
Kipnis~\cite{FerrariKipnis} it can be proved that $X_0(t)/t$ converges in
distribution to a random variable uniform in $[-\rho,\rho]$, where
hereafter we denote $\rho=2p-1$. Note that the particles in the exclusion
process try to perform a random walk with drift $\rho$ (and they cannot go
faster than that), that explains why the support of the limiting random
variable is changed. In fact, in many ways the ASEP behaves similarly to
the TASEP slowed down by a factor of $\rho$.
\begin{conj}\label{C:ASEP_lln}
In the ASEP, $\lim_{t\to\infty} X_0(t)/t$ exists a.s. (and the limit is
uniform on $[-\rho,\rho]$).
\end{conj}

By the discussion preceding Corollary~\ref{C:speed} this is equivalent to
the \mbox{following}:

\begin{conj}
The ASEP speed process measure $\mu^{\mathrm{ASEP}}$ is well defined and
translation invariant with each $U_n$ uniform on $[-\rho,\rho]$.
\end{conj}

In order for statements about the \textit{ASEP} speed process to make
sense we
must assume this conjecture, and therefore some of our theorems are
conditional on Conjecture~\ref{C:ASEP_lln}. It should be noted that
with minor
modifications our results also hold assuming a weaker assumption,
namely a
joint limit in distribution of the speeds $\{X_n(t)/t\}_{n\in\Z}$. In that
case, the speed process measure is still defined, even though the particles
may not actually have an asymptotic speed.

As noted there, Theorem~\ref{T:symmetry} holds also for the ASEP, with no
additional condition. Theorem~\ref{T:stationary} becomes conditional:
\begin{theorem}\label{T:stationary_ASEP}
Assume Conjecture~\ref{C:ASEP_lln} holds. Then $\mu^{\mathrm{ASEP}}$ is a stationary
measure for the ASEP: the unique ergodic stationary measure which has
marginals uniform on $[-\rho,\rho]$.
\end{theorem}

As in the case of the TASEP, this can be interpreted as follows: if an ASEP
is started with initial configuration in $[-\rho,\rho]^\Z$ with
distribution $\mu^{\mathrm{ASEP}}$, then at any time $t>0$\vadjust{\goodbreak} the distribution of the
process is also given by $\mu^{\mathrm{ASEP}}$. Note that both the dynamics and the
measure $\mu^{\mathrm{ASEP}}$ depend implicitly on the asymmetry parameter~$p$.

A useful tool in studying the speed process is the understanding of the
stationary measures of the $N+1$ type TASEP in terms of a multi-line
process described below, developed by Angel~\cite{Angel} and Ferrari and
Martin~\cite{FerrariMartin}. There is no known analogue for these results
that describes the stationary measure of the multi-type ASEP. Thus we need
to use other (and weaker) techniques to extract information about the
marginals of the ASEP speed process. This explains the contrast in the
level of detail between the following results and the corresponding
theorems above about the TASEP.
\begin{theorem}\label{T:ASEP-swap}
We have the following limit:
\[
\lim_{t\to\infty} \P\bigl(X_0(t) < X_1(t)\bigr) = \frac{2-p}{3}.
\]
\end{theorem}

Theorem 2.3 of~\cite{FGM} proves that the probability that particles $0$
and $1$ interact at least once (i.e., one of them tries to jump onto the
other) is $\frac{1+p}{3p}$. In the next section we will show that this is
equivalent to the just stated theorem.

Our next theorem provides information about the joint distribution of
$\{U_0, U_1\}$, assuming Conjecture~\ref{C:ASEP_lln} holds.
\begin{theorem}\label{T:ASEP-joint}
Assume Conjecture~\ref{C:ASEP_lln} holds. Let the measure $\muI$ on
$[-\rho,\rho]^2$ be the marginal of $\{U_0,U_1\}$ under $\mu^{\mathrm{ASEP}}$.
Denote by $\muII$ the reflection of $\muI$ about the line $x=y$. Then on
$\{(x,y)\dvtx-\rho\le x<y\le\rho\}$ we have\looseness=-1
\[
p \cdot\muI- \bar{p} \cdot\muII=
\frac{y-x}{4\rho^2} \,dx \,dy.
\]\looseness=0
\end{theorem}

We finish this section with a statement concerning the case $U_0=U_1$.
Consider the total amount $J_{i,j}$ of time that particles $i$ and $j$
spend next to each other, that is, $J_{i,j} = \int_0^\infty
\ind(|X_i(t)-X_j(t)|=1) \,dt$.
\begin{theorem}\label{T:together_forever}
In the TASEP, $J_{0,1}=\infty$ if and only if $U_0=U_1$. If
Conjecture~\ref{C:ASEP_lln} holds, then the same holds for the ASEP.
\end{theorem}

In the TASEP $J_{0,1}=\infty$ implies that there is at least one
interaction between 0 and 1 which means that they are a.s. swapped. (See
the next section for a more detailed discussion.) Thus if $U_0=U_1$, then
eventually $X_0(t)>X_1(t)$. In fact, this holds for any two particles in
the same convoy: in Lemma~\ref{L:diagswap} we will prove that in the TASEP,
particle 0 will eventually overtake all the particles in its convoy with
positive index.

\subsection{Overview of the paper}

The rest of the paper is organized as follows. Section~\ref{S:background}
provides some of the background: constructions of the processes and the
multi-line description of\vadjust{\goodbreak} the stationary measure for the multi-type TASEP.
Section~\ref{S:symmetry} includes the proof of the symmetry property
(Theorem~\ref{T:symmetry}) and Section~\ref{S:stationary} proves the
stationarity of
the speed process (Theorems~\ref{T:stationary} and
\ref{T:stationary_ASEP}). Sections~\ref{S:joint} and~\ref{S:mult} include
the results about various finite-dimensional marginals of the TASEP speed
process. Section~\ref{S:convoys} deals with the proof of Theorem \ref
{T:convoys}.
Finally, in Section~\ref{S:ASEP} we prove our results on the ASEP speed process.

\section{Preliminaries} \label{S:background}

\subsection{Construction of the process}

There are several formal constructions of the TASEP and ASEP. The one that
best suits our needs seems to be Harris's approach~\cite{Liggett}. We
include the construction since there are several variations and the exact
details are used in some of our proofs. The process is a function $Y$
defined on $\Z\times\R^+$. $Y_k(t)$ will denote the class of the particle
at position $k$ at time $t$. The configuration at time $t$ is $Y(t) =
\{Y_k(t)\}_{k\in\Z}$. The classes\vspace*{1pt} of particles will be real numbers, hence
the configuration at any given time is in $\R^\Z$. Setting $t=0$ gives the
initial configuration $Y(0)$.

We define the transposition operator $\tau_n$, acting on $\R^\Z$ by
exchanging $Y_n$ and $Y_{n+1}$, while keeping all other classes equal.
Using this we can alternately describe the sorting operator $\sigma_n$ by
\[
\sigma_n Y = \cases{
\tau_n Y, &\quad$Y_n<Y_{n+1}$, \cr
Y, &\quad otherwise.}
\]
Thus $\sigma_n$ has the effect of sorting $Y_n,Y_{n+1}$ in decreasing
order, keeping other classes the same.

The TASEP is defined using the initial configuration and the location of
``jump'' points. The probability space contains a standard Poisson process
on $\Z\times\R^+$, that is, a collection of independent standard Poisson
processes on $\R^+$, denoted $T_n$. If $(n,t)$ is a point of $T_n$,
then at
time $t$ the values of $Y_{n}(t^-)$ and $Y_{n+1}(t^-)$ may be switched. In
the TASEP they are sorted, that is, $Y(t) = Y(t^-)\cdot\sigma_n$. This
can be
described as applying each of the operators $\sigma_n$ at rate 1
independently. A simple percolation argument shows that this dynamic is
a.s. well defined. (For any fixed $t>0$ there are a.s. infinitely many
integers $n$ so that there are no Poisson points on $\{n\}\times[0,t]$
which means that to define the process up to time $t$ it suffices to
consider finite lattices.)

\subsubsection*{The ASEP}
Defining the partially asymmetric exclusion process requires additional
randomness. Given the parameter $p\in(1/2,1]$, we attach to each point
$(n,t)$ in the Poisson process an independent Bernoulli random variable
$X_{n,t}$ with $\P(X_{n,t}=1) = p$. We can now
define the probabilistic sorting operator $\rho_n$ as follows:
\[
\rho_n Y = \cases{
\sigma_n Y, &\quad if $X_{n,t}=1$, \cr
\sigma^*_n Y, &\quad if $X_{n,t}=0$.}
\]
Thus with probability $p$ the smaller classed particle is moved to the
right position and with probability $1-p$ it is moved to the left position.
When such an event happens we say that $Y_n(t)$ and $Y_{n+1}(t)$ have an
\textit{interaction} (regardless of whether they were actually
swapped). Note
that if particles $i,j$ interact in this way, then their order after the
swap is independent of the order before the swap. The key observation is
that after $i<j$ interact in this way at least once, $i$ has probability
$p$ of being to the right of $j$, and this is unchanged by further
interactions. Moreover, if we condition on $J_{i,j}(t)= \int_0^\infty
\ind(|X_i(s)-X_j(s)|=1) \,ds$ (the total time $i,j$ spend next to each other
until time $t$), then
%
%
\begin{equation}\label{eq:J}\quad
\P\bigl(X_i(t)<X_j(t)| J_{i,j}(t)\bigr)=e^{-J_{i,j}(t)} + \bar p
\bigl(1-e^{-J_{i,j}(t)}\bigr)=\bar p + p e^{-J_{i,j}(t)},
\end{equation}
where the expression on the right is just the probability that there were
no interaction between $i$ and $j$ until time $t$ plus the probability that
there was some interaction, and at time $t$ particle $i$ is to the left of
$j$. One of the consequences of (\ref{eq:J}) is that
%
%
\begin{equation}\label{eq:togetherJ}
\lim_{t\to\infty} \P\bigl(X_i(t)<X_j(t)\bigr)=\bar p + p \E e^{-J_{i,j}}.
\end{equation}
Thus Theorem~\ref{T:ASEP-swap} implies $\bar p+p \E e^{-J}=\frac{2-p}3$ which
gives $1-\E e^{-J_{0,1}} = \frac{1+p}{3p}$. But $1 - \E e^{-J_{0,1}}$ is
exactly the probability that there is at least one interaction between $0$
and 1 which shows why Theorem 2.3 of~\cite{FGM} and our Theorem \ref
{T:ASEP-swap}
are equivalent.

In the TASEP case if there is an interaction between $i<j$, then
$X_i(t)>X_j(t)$ after that. Thus in that case from (\ref{eq:togetherJ}) we
get
\[
\P\bigl(\mbox{eventually } X_i(t)>X_j(t) | J_{i,j} = \infty\bigr) = 1,
\]
which explains the remark after Theorem~\ref{T:together_forever}.

There is an alternate construction for the ASEP, which will be used in
Section~\ref{S:symmetry}. Consider a Poisson process with lower intensity
$p$ on $\Z\times\R^+$, but whenever it has a point $(n,t)$ we apply at time
$t$ the operator $\pi_n$ rather then $\rho_n$, where $\pi_n$ is defined by
\[
\pi_n Y = \cases{
\tau_n Y, &\quad$Y_n<Y_{n+1}$, \cr
\tau_n Y, &\quad$Y_n>Y_{n+1}$ with prob. $q=(1-p)/p$, \cr
Y, &\quad$Y_n>Y_{n+1}$ with prob. $\bar q = (2p-1)/p$.}
\]
Thus if the pair is in increasing order it is always swapped, while if it
is in decreasing order it is swapped only with probability $q$. It is easy
to see that every possible swap occurs at the same rate in the two
constructions; hence the resulting processes have the same generator.

\subsection{Stationary measures for the multi-type TASEP}
\label{S:stat_tasep}

The following theorem can be proved by standard coupling methods (see,
e.g.,
\cite{Liggettcouple} where the same theorem is proved for the 2-type
TASEP).\vadjust{\goodbreak}
\begin{theorem}\label{T:multistat}
Fix every $0 \leq\lambda_1, \ldots, \lambda_N \leq1$ with $\sum
\lambda_i = 1$. There is a unique ergodic stationary distribution
$\nu_\lambda$ for the $N$-type TASEP with $\P(Y_0=k)=\lambda_k$. The
measures $\nu_\lambda$ are the extremal stationary translation invariant
measures. They are the only stationary translation invariant measures
with the property that for each $k$, the distribution of $\{ \ind
[Y_n\leq
k] \}_{n\in\Z}$ is product Bernoulli measure with density $\sum_{j\leq k}
\lambda_j$.
\end{theorem}

For the ordinary TASEP (with particles and holes) this stationary
distribution is just the product Bernoulli with a fixed density. If we have
an $(N+1)$-type TASEP then the structure of the stationary distribution is
more complicated. The first description of $\nu_\lambda$ for $N=2$ was
given by the matrix method~\cite{DJLS}. Reference~\cite{FFK} gave probabilistic
interpretations and proofs of the measure and its properties. Recently
combinatorial descriptions of $\nu_\lambda$ have appeared as well. The
$(2+1)$-type TASEP was treated by Angel~\cite{Angel} (see also
Duchi--Schaeffer~\cite{DuchiSchaeffer}). These results were extended
for all
$N$ by Ferrari and Martin~\cite{FerrariMartin}. They give an elegant
construction of $\nu_\lambda$ using systems of queues.

We will now briefly describe the $N$-line description of $\nu_\lambda$ for
the $(N+1)$-type TASEP. The two-line case suffices for most of our results,
with the exception of the results of Section~\ref{S:mult}. For a more
detailed description and proofs see~\cite{FerrariMartin}.

From here on we shall fix the parameters $\lambda_1,\ldots,\lambda_{N+1}$.
Consider $N$ independent Bernoulli processes on $\Z$ denoted $\bB_1, \bB_2,
\ldots, \bB_N$ where $\bB_k$ has parameter $\sum_{i\leq k} \lambda_i$ (these
are the lines). From these lines we construct a system of $N-1$ coupled
queues. The lines give the service time of the queues, and the departures
from each queue are the arrivals to the next queue.

It is important to observe that the time for the queues goes from right to
left, that is, $\bB_i(n)$ is followed by $\bB_i(n-1)$ and so on. The resulting
system of queues is positively recurrent, so it can be defined starting at
$\infty$ and going over the lines toward $-\infty$.

The $i$th queue will consist of the particles that departed from the
$i$th line and are waiting for a service in $\bB_{i+1}$. This queue will
consist of particles of classes $\{1,\ldots,i\}$. When a service is
available in $\bB_{i+1}$ the lowest classed particle in the $i$th queue is
served and departs (to the next queue). If the queue is empty then a
particle of class $i+1$ is said to depart the queue. The departure process
of each queue (i.e., the times and sequence of classes of departing
customers) is the arrival process for the next queue.

It is convenient to think of an additional queue with $\bB_1$ as its
service times. This queue has no arrivals (so it is always empty). The
unused services introduce first class particles, which join the second
queue whenever there is a service in~$\bB_1$. These operations are
evaluated for each $n$ from line $1$ to line $N$ in order. Let $Q_{i,j}(n)$
be the number of particles of type $j$ in the $i$th queue after column $n$
of the multi-line process has been used.

Note that each queue has a higher rate of service than of arrivals, so the
queues sizes are tight, and the state with all queues empty is positively
recurrent. In practice, the $i$th queue has $i$ types of particles in it,
so the whole system of queues is described by $\frac{N(N-1)}{2}$
nonnegative integers.
\begin{theorem}[(Ferrari--Martin)]
$\nu_\lambda$ is the distribution of the departure process of $\bB_N$,
with class $N+1$ (or empty sites) at those $n$ when there is no service.
\end{theorem}

As an example, and to clarify the graphic representation we use later,
consider the following segment of a configuration of the three-line
process for
\mbox{$n=\{1,2,3,4\}$}. Suppose both queues are empty at time 5. (This is denoted
by the $\varnothing,\varnothing$ exponent.) Here, $\C$ denotes a $0$ in the
corresponding line, and $\D$ a 1. Later, in cases where we do not care
about a specific value we may use $\B$ to denote that
\[
\config{4}{{\C}&\D&\D&{\C}\\\D&\C&\C&\D\\\D&\D&\C&\D}^{\varnothing
,\varnothing}.
\]
At time $4$, reading the rightmost column from top to bottom, there is no
service in $\bB_1$, so no first class particle joins the second queue,
which therefore remains empty. There is a service in $\bB_2$, and no
particles in the first queue, so a second class particle joins the second
queue. There is service in $\bB_3$, so the second class particle departs
immediately. Thus at time $4$ the queue states are $(\varnothing
,\varnothing)$.

At time $3$ a first class particle arrives to the first queue, and
stays there since there is no service in the second queue. There is
no further service in column $3$, so the state at time $3$ is
$(\{1\},\varnothing)$. There is no departure, which is denoted by a
$4$ (or hole). At time $2$ another first class particle arrives, and
there is no particle in the second queue so the service in $\bB_3$
gives rise to a third class particle departing. The states are now
$(\{1,1\},\varnothing)$. Finally, at time $1$ a first class particle
is served at both $\bB_2$ and $\bB_3$, departing and leaving queue
states $(\{1\},\varnothing)$. The resulting segment of $\nu_\lambda$
is $(1,3,4,2)$.

\section{Symmetry} \label{S:symmetry}

Recall the operators $\pi_n$ defined above. These act randomly on
configurations, and the ASEP can be defined by applying each of the Markov
operators $\pi_n$ at rate $p$.

Formally, $\pi_i$ is defined as acting on $\cM(S_\infty)$: probability
measures on $S_\infty$. Given a measure $\nu$ on $S_\infty$, we let
$\pi_n\nu$ be the distribution of $\pi_n$ applied to a sample from $\nu$.
Since $\tau_i$ and $\sigma_i$ also act naturally on the measures (in the
same way), one finds the operator relation
\[
\pi_i = q\tau_i + {\bar q} \sigma_i.
\]
Note that $p=1$ gives $q=0$ so in that case $\pi_i=\sigma_i$.
In the case $p=1/2$ we get $q=1$ and $\pi_i=\tau_i$, so the process reduces
to a symmetric random walk on $S_\infty$.\vadjust{\goodbreak}

The crucial observation leading to Theorem~\ref{T:symmetry} is the following
lemma.
\begin{lemma}\label{L:reverse}
Fix any $p\ge1/2$, and sequence $i_1,\ldots,i_n$. Then
%
%
\begin{equation}
\label{eq:reverse}
\pi_{i_n} \cdots\pi_{i_1} \cdot\id
\eqd( \pi_{i_1} \cdots\pi_{i_n} \cdot\id)^{-1}.
\end{equation}
\end{lemma}

That is, applying a sequence of $\pi_i$'s in the reverse order to the
identity leads to the inverse permutation. This is trivially true when
$p=1/2$ and $\pi=\tau$, but requires proof for other $p$. When
$p\in\{1/2,1\}$ the operator is deterministic and this distributional
identity is an equality of permutations.
\begin{pf*}{Proof of Theorem~\ref{T:symmetry}}
The theorem follows from Lemma~\ref{L:reverse} since at any finite
time at
each $i$ there is positive probability ($e^{-t}$) that no swap has
occurred. Each such $i$ separates $\Z$ into two parts with independent
behavior, so the state of the process is a product of finite, mutually
commuting permutations. The distribution of the sequence of applied
operators between such inactive locations is symmetric in time.
\end{pf*}

We now prove Lemma~\ref{L:reverse}. In the case of the TASEP, Lemma \ref
{L:reverse}
and Theorem~\ref{T:symmetry} were first proved in~\cite{DPRW}. To prove the
lemma in the general case, we start with the following facts about the
transposition operators. The identities are readily verified, and the last
claim is known as Matsumoto's lemma (see, e.g.,
\cite{coxetergroups}, Theorem 3.3.1).
\begin{fact}
The operators $\tau_i$ satisfy the relations
%
%
\begin{eqnarray}
\label{rel1}
\tau_i^2 &=& I, \\
\label{rel2}
\tau_i\tau_j &=& \tau_j\tau_i \qquad\mbox{for }
|i-j|>1, \\
\label{eq:rel3}
\tau_i\tau_{i+1}\tau_i &=& \tau_{i+1}\tau_i\tau_{i+1},
\end{eqnarray}
where $I$ denotes the identity operator. With these relations the
operators $\{\tau_i\}$ generate the symmetric group. Furthermore, it is
possible to pass between any two minimal words of the same permutation
(i.e., words of minimal length representing that permutation) using only
relations \eqref{rel2}, \eqref{eq:rel3}.
\end{fact}

The $\pi$'s satisfy similar relations:
\begin{lemma}\label{lemma3.3}
The operators $\{\pi_i\}$ satisfy the relations
%
%
\begin{eqnarray}
\label{eq:rel1p}
\pi_i^2 &=& q I + \bar q \pi_i, \\
\label{eq:rel2p}
\pi_i\pi_j &=& \pi_j\pi_i\qquad \mbox{for }
|i-j|>1, \\
\label{eq:rel3p}
\pi_i\pi_{i+1}\pi_i &=& \pi_{i+1}\pi_i\pi_{i+1}.
\end{eqnarray}
\end{lemma}

Note that only the first relation differs from the corresponding relation
for $\tau$. When $p=1/2$ these reduce to the relations for $\tau$. In the
case $p=1$ the first relation becomes $\sigma_i^2=\sigma_i$.\vadjust{\goodbreak} In that case,
the only nontrivial relation is~\eqref{eq:rel3p} which is true since both
sides have the effect of sorting the three terms involved in decreasing
order.
\begin{pf*}{Proof of Lemma~\ref{lemma3.3}}
Equation \eqref{eq:rel1p} is easy to check, and \eqref{eq:rel2p} is trivial.
For~\eqref{eq:rel3p}, using $\pi= q\tau+ \bar q \sigma$ and expanding, we
need to show that
\begin{eqnarray*}
&&q^3 (\tau_i\tau_{i+1}\tau_i) +
q^2\bar q (\tau_i\tau_{i+1}\sigma_i + \tau_i\sigma_{i+1}\tau_i +
\sigma_i\tau_{i+1}\tau_i) \\
&&\qquad{} + q \bar q^2 (\tau_i\sigma_{i+1}\sigma_i + \sigma_i\tau
_{i+1}\sigma_i +
\sigma_i\sigma_{i+1}\tau_i) +
\bar q^3 (\sigma_i\sigma_{i+1}\sigma_i)
\end{eqnarray*}
is unchanged by exchanging $i$ and $i+1$. It is easy to verify that
\begin{eqnarray*}
\tau_i\tau_{i+1}\tau_i &=& \tau_{i+1}\tau_i\tau_{i+1},\qquad
\tau_i\tau_{i+1}\sigma_i = \sigma_{i+1}\tau_i\tau_{i+1} ,\qquad
\tau_i\sigma_{i+1}\tau_i = \tau_{i+1}\sigma_i\tau_{i+1}, \\
\sigma_i\sigma_{i+1}\sigma_i &=&
\sigma_{i+1}\sigma_i\sigma_{i+1},\qquad
\sigma_i\tau_{i+1}\tau_i = \tau_{i+1}\tau_i\sigma_{i+1},
\end{eqnarray*}
so it remains to show
\[
\tau_i\sigma_{i+1}\sigma_i + \sigma_i\tau_{i+1}\sigma_i +
\sigma_i\sigma_{i+1}\tau_i =
\tau_{i+1}\sigma_i\sigma_{i+1} + \sigma_{i+1}\tau_i\sigma_{i+1} +
\sigma_{i+1}\sigma_i\tau_{i+1}.
\]
We may assume $i=0$. Since only the relative order of
$\eta_0,\eta_1,\eta_2$ matters, we may assume these are $\{0,1,2\}$ in
some order. Applying these operators to the $6$ possible orders  gives Table~\ref{table1}.
%
\begin{table}
\tablewidth=240pt
\caption{See proof of Lemma \protect\ref{lemma3.3}}\label{table1}
\begin{tabular*}{\tablewidth}{@{\extracolsep{\fill}}lcccccc@{}}
\hline
$\bolds{\eta}$ & \textbf{012} & \textbf{021} & \textbf{102} & \textbf{120} & \textbf{201} & \textbf{210} \\
\hline
$\tau_0\sigma_1\sigma_0 \cdot\eta$ & 210 & 120 & 210 & 120 & 120 &
120 \\
$\sigma_0\sigma_1\tau_0 \cdot\eta$ & 210 & 210 & 201 & 210 & 201 &
210 \\
$\sigma_0\tau_1\sigma_0 \cdot\eta$ & 210 & 210 & 210 & 201 & 210 &
201 \\
[6pt]
$\tau_1\sigma_0\sigma_1 \cdot\eta$ & 210 & 210 & 201 & 201 & 201 &
201 \\
$\sigma_1\sigma_0\tau_1 \cdot\eta$ & 210 & 120 & 210 & 120 & 210 &
210 \\
$\sigma_1\tau_0\sigma_1 \cdot\eta$ & 210 & 210 & 210 & 210 & 120 &
120 \\
\hline
\end{tabular*}
\end{table}
In each column, the entries in the top half are a permutation of the
entries in the bottom half, so adding the first three operators
gives the same result as adding the last three.
\end{pf*}
\begin{pf*}{Proof of Lemma~\ref{L:reverse}}
Given $(i_1,\ldots,i_n)$, let $X = \tau_{i_1} \cdots\tau_{i_n}$. If this
is a minimal (w.r.t. length) word for $X$ in $S_\infty$, then $\pi_{i_n}
\cdots\pi_{i_1} \cdot\id= X$ with probability 1. In this case, the
reverse word is minimal for $X^{-1}$, so the claim holds.

The proof proceeds by induction on $n$. Take some sequence
$(i_1,\ldots,i_n)$. If the representation $\tau_{i_1} \cdots\tau
_{i_n}$ is minimal,
then the claimed identity holds. Otherwise, let $k$ be maximal
such that $X = \tau_{i_1} \cdots\tau_{i_k}$ is a minimal representation.
By maximality of $k$ we see that $Y = X \tau_{i_{k+1}}$ has a shorter
representation, so there is a representation $Y = \tau_{j_1} \cdots
\tau_{j_{k-1}}$. (The length is $k-1$ and not $k$ since its parity is
opposite that of~$X$.) Thus $X = \tau_{j_1}\cdots\tau_{j_{k-1}}
\tau_{i_{k+1}}$ is another minimal representation of $X$.\vadjust{\goodbreak}

Starting with $\pi_{i_1} \cdots\pi_{i_n}$, we can repeatedly apply
relations \eqref{eq:rel2p} and \eqref{eq:rel3p} to the first $k$ terms in
the product, to get
\[
\pi_{i_1} \cdots\pi_{i_n} = \pi_{j_1} \cdots\pi_{j_{k-1}}
\pi_{i_{k+1}}^2 \cdots\pi_{i_n}.
\]
Here $i_{k+1}^2$ appears twice since it is both the last term in the
alternate representation of $X$ and the first in the remainder of the
sequence. Relation \eqref{eq:rel1p} now gives
%
%
\begin{equation}\label{eq:short1}\quad
\pi_{i_1} \cdots\pi_{i_n} = q (\pi_{j_1} \cdots\pi_{j_{k-1}}
\pi_{i_{k+2}} \cdots\pi_{i_n}) + \bar{q} (\pi_{j_1} \cdots
\pi_{j_{k-1}} \pi_{i_{k+1}} \cdots\pi_{i_n}).
\end{equation}
Similarly, working with the reverse sequence,
%
%
\begin{equation}\label{eq:short2}\quad
\pi_{i_n} \cdots\pi_{i_1} = q (\pi_{i_n} \cdots\pi_{i_{k+2}}
\pi_{j_{k-1}} \cdots\pi_{j_1}) + \bar{q} (\pi_{i_n} \cdots
\pi_{i_{k+1}} \pi_{j_{k-1}} \cdots\pi_{j_1}).
\end{equation}
Applying \eqref{eq:short2} and \eqref{eq:short1} to to $\id$, and using
the induction hypothesis for the shorter sequences $(j_1, \ldots,
j_{k-1}, i_{k+1},\ldots, i_n)$ and $(j_1, \ldots, j_{k-1},
i_{k+2},\ldots,
i_n)$ completes the proof.
\end{pf*}

\textit{Note}: the proof actually shows that any word in the $\pi$'s
can be reduced
(as an operator) to some convex combination of words corresponding to
minimal words.
\begin{coro}\label{C:lim}
Consider the infinite type TASEP with initial condition $Y_n(0)=n$. Then
$\{\frac{Y_n(t)}{t}\}_{n\in\Z}$ converges weakly to $\mu$ as
$t\to\infty$.
\end{coro}
\begin{pf}
For any $t$ this process has the same law as $\{\frac{X_n(t)}{t}\}
_{n\in
\Z}$, which converges a.s. to a process with law $\mu$.
\end{pf}

\section{Stationarity} \label{S:stationary}

We will give two different proofs of the stationarity of the distribution
of the speed process. The first is specific to the TASEP, and is
reminiscent of coupling from the past. It uses the Harris construction
directly. The second proof is based on the symmetry between $\{X_n(t)\}$
and $\{Y_n(t)\}$ (or more specifically Corollary~\ref{C:lim}). The second
proof holds also for the ASEP, word by word, under the assumption that
Corollary~\ref{C:lim} is true for the ASEP (which is weaker then
Conjecture~\ref{C:ASEP_lln}).

\subsection{Coupling proof}

\begin{lemma}
Consider two TASEPs $Y, Y'$ defined via the Harris construction as the
function of the same Poisson process on $\Z\times\R^+$. We set the
initial conditions as $Y_n(0)=n$ and $Y'(0)= \sigma_0 Y(0)$ (i.e.,
particles 0 and 1 are switched initially in~$Y'$). Let
$\{U_n\}=\{\lim_{t\to\infty} X_n(t)/t\}$ denote the speed process of
$Y$, and $\{U'_n\}=\{\lim_{t\to\infty} X'_n(t)/t\}$ denote the speed
process of $Y'$. Then $U' = \sigma_0 U$.
\end{lemma}
\begin{pf}
All particles other than $\{0,1\}$ are either larger or smaller than both
$0$ and $1$, so any swaps involving a particle other than $\{0,1\}$ will
occur or not\vadjust{\goodbreak} occur equally in $Y$ and $Y'$. It follows that for any
$i\notin\{0,1\}$ we have $X_i(t)=X'_i(t)$ and hence $U_i=U'_i$. Similarly,
since 0 and 1 must fill the only vacant trajectories,
$\{U_0,U_1\}=\{U'_0,U'_1\}$ as an unordered pair.

In $Y'$ particle $0$ is always to the right of particle $1$, so
$U'_0=\max\{U_0,U_1\}$ and $U'_1=\min\{U_0,U_1\}$, completing the proof.
\end{pf}
\begin{pf*}{Proof of Theorem~\ref{T:stationary} using coupling}
Consider a Poisson process on $\Z\times\R$. Half of the process, namely
the restriction to $\Z\times\R^+$ is used in the Harris construction of
the TASEP. Similarly, for any $s\in\R$ we can translate the Poisson
process by $s$ [i.e., take all points of the form $(n,t+s)$ where $(n,t)$
is in the original process], and take the restriction to $\Z\times\R^+$,
which can be used in the Harris construction to get a different (though
highly dependent) instance of the TASEP.

Let $U_n(s)$ be the speed process resulting from the Harris construction
using the translated Poisson process. Clearly for every $s$, $U(s)$ has
the same law $\mu$, so we are done if we show that $U_n(s)$ evolves as a
TASEP (with time parameter $s$). Consider the effect of an infinitesimal
positive shift $s$. The shift adds new $\sigma$ operations, to be applied
before the original sequence of operations. These are added at rate~1 at
each location. By the previous lemma, the effect on the resulting speeds
of applying $\sigma_n$ before using the same Poisson process is to apply
$\sigma_n$ to the speeds, which is exactly what we need.\looseness=-1
\end{pf*}

It is interesting to note that in the Poisson process $\Z\times\R$, the
part on $\Z\times\R^+$ is used to determine the ``initial'' speed process
$U(0)$, and the restriction to $\Z\times\R^-$ is used exactly as in the
Harris construction to generate the TASEP dynamics of $U(s)$.

\subsection{Symmetry based proof}

\mbox{}

\begin{pf*}{Proof of Theorems~\ref{T:stationary} and
\ref{T:stationary_ASEP} using symmetry}
We write the proof for $\mu$, but it holds verbatim for $\mu^{\mathrm{ASEP}}$
under Conjecture~\ref{C:ASEP_lln}.

Informally, we argue as follows. Fix $s$ and let $t\to\infty$. Both
$\frac{X_t}{t}$ and $\frac{X_{t+s}}{t+s}$ converge a.s. to a sample of
$\mu$. By Theorem~\ref{T:symmetry} these have the same law as $\frac{Y_t}{t}$
and $\frac{Y_{t+s}}{t+s}$, so for large $t$ both of these have law close
to $\mu$. However, the result of letting $\frac{Y_t}{t}$ evolve for an
additional $s$ time is $\frac{Y_{t+s}}{t}$, which is close to
$\frac{Y_{t+s}}{t+s}$.

Let $P_s$ be the evolution operator for the Markov process corresponding
to the generator $L$ on $\R^{\Z}$ [see (\ref{eq:gen1})]. To
prove\vspace*{1pt}
stationarity it is enough to show that for every $0<s$ and every bounded
continuous local function $f\dvtx\R^{\Z}\to\R$ we have
%
%
\begin{equation}\label{eq:statmu}
\int P_s f(\eta) \,d\mu(\eta) = \int f(\eta) \,d\mu(\eta).
\end{equation}
Consider the process $\{Y_n(t)\}_{n\in\Z}$ started from $Y_n(0)=n$ and
denote the distribution of $\{ \frac{Y_n(t)}{t}\}_n$ by
$\nu_t$. By Corollary~\ref{C:lim}\vadjust{\goodbreak} the weak limit of $\nu_t$ is $\mu$
which means that for every local bounded continuous function
$f\dvtx\R^{\Z}\to\R$ we have
\[
\int f(\eta) \,d\nu_t(\eta) \xrightarrow{t\to\infty}
\int f(\eta) \,d\mu(\eta).
\]
For any fixed $s$
\[
\int P_s f(\eta) \,d\nu_t(\eta) \xrightarrow{t\to\infty}
\int P_s f(\eta) \,d\mu(\eta).
\]
But $\int P_s f(\eta) \,d \nu_t(\eta) = \int f(\frac{t+s}{t} \eta
)\, d \nu_{t+s}(\eta)$ which (for any fixed $s$, as $t\to\infty$)
converges to $\int f(\eta) \,d\mu(\eta)$. Now \eqref{eq:statmu} and the
theorem follow.
\end{pf*}

\section{Basic properties of stationary distributions}

In this section we present a medley of simple results concerning the
(T)ASEP and its stationary distributions. These are only weakly related to
each other, and are collected here for convenience.
\begin{prop}\label{P:erg}
$\mu$ is ergodic for the shift. Under Conjecture~\ref{C:ASEP_lln}, so is
$\mu^{\mathrm{ASEP}}$.
\end{prop}
\begin{pf}
Consider the setup of Corollary~\ref{C:speed} and use the Harris
construction with independent standard Poisson processes $T_n$ on the
interval $[0,\infty)$ to define $y_n(t)$ and the variables $X_n(t)$. Then
the limit process $\{U_n\}_{n\in\Z}$ is measurable with respect to the
$\sigma$-algebra $\mathcal{F}$ generated by the i.i.d. processes $T_n$
($n\in\Z$). Since $\mathcal{F}$ is generated by i.i.d. processes any
translation invariant event in $\mathcal{F}$ has to be trivial. But then
the same thing must be true for any translation invariant event in the
$\sigma$-algebra generated by $\{U_n\}_{n\in\Z}$ as this is a
sub-$\sigma$-algebra of $\mathcal{F}$.
\end{pf}

There are three possible ``reflections'' for the ASEP. One may reverse the
direction of space, so that (low classed) particles flow to the left and
not right; one can consider the time reversal of the dynamics, and one can
reverse the order of classes (or keep the same generator but replace class
$k$ with $-k$, or $N+1-k$, etc.). It is easy to see that reversal of both
space and class order preserves the original dynamics. This is called the
space-class symmetry of the TASEP/ASEP.

The following proposition is the \textit{space-class} symmetry of the speed
process, and follows directly from the corresponding symmetry of the ASEP
process.
\begin{prop}\label{P:spacesym}
For the TASEP $\{U_n\}_{n\in\Z} \eqd\{-U_{-n}\}_{n\in\Z}$. This also
holds for the ASEP, assuming Conjecture~\ref{C:ASEP_lln} holds.
\end{prop}

The following observation and its corollary provide an important connection
between the distribution of the speed process and the stationary measures
of multi-type ASEP.\vadjust{\goodbreak} These connections will be used to extract information
on the joint distribution of the speeds of several particles in
Sections~\ref{S:joint} and~\ref{S:mult}.
\begin{lemma}\label{L:projection}
Let $\{\eta_n(t)\}_{n,t}$ be an ASEP, and let $F\dvtx\R\rightarrow\R$
be a
nondecreasing function. Then $\{F(\eta_n(t))\}_{n,t}$ is also an ASEP
(with the same asymmetry parameter).
\end{lemma}
\begin{pf}
The ASEP is defined as applying to $\eta(t)$ each of the operators
$\pi_n$ independently at rate $1$. Applying a nondecreasing function to
each coordinate commutes with every $\pi_i$, hence
$\{F(\eta_n(t))\}_{n,t}$ is just the ASEP with initial configuration
$\{F(\eta_n(0)\}_n$.
\end{pf}
\begin{coro}\label{C:projection_to_N}
If $F\dvtx[-1,1]\rightarrow\{1,\ldots,N\}$ is nondecreasing, then for the
TASEP the distribution of $\{F(U_n)\}$ is the unique ergodic stationary
measure of the multi-type TASEP with types $\{1,\ldots,N\}$ and densities
$\lambda_i = \frac{1}{2}\Leb(F^{-1}(i))$.

This also holds for the ASEP (and its corresponding multi-type stationary
measure) under Conjecture~\ref{C:ASEP_lln}.
\end{coro}
\begin{pf}
Let $\mu_F$ denote the distribution of $\{F(U_n)\}$. Since $\mu$ is
ergodic, so is~$\mu_F$. The marginals are as claimed since each $U_n$ is
uniform on $[-1,1]$.

To prove that $\mu_F$ is stationary, start a TASEP $Y_n(t)$ with initial
configuration $Y_n(0)=U_n$. By Lemma~\ref{L:projection} $\{F(Y_n(t))\}_{n,t}$
is a $N$-type TASEP. Since $\mu$ is stationary, $Y(t)$ also has law
$\mu$, and so $\{F(Y_n(t))\}_{n,t} \eqd\{F(Y_n(0))\}_{n,t}$, hence
$\mu_F$ is also stationary.

The result for the ASEP follows the same way.
\end{pf}

The next proposition shows that a TASEP started with uniform i.i.d.
classes must converge to the speed process. In particular, even though
classes in the i.i.d. initial distribution are a.s. all different, the
process converges to the speed process which has infinite convoys of
particles with the same class (see Section~\ref{S:convoys}). Thus the TASEP
dynamics has the effect of aggregating particles with increasingly closer
speeds next to each other.
\begin{prop}
Consider a TASEP where $Y_n(0)$ are i.i.d. uniform on $[-1,1]$. Then
$\{Y_n(t)\}_{n\in\Z}$ converges weakly to $\mu$. The same holds for the
ASEP under Conjecture~\ref{C:ASEP_lln}
\end{prop}
\begin{pf}
Let $\nu_t$ be the distribution of $Y(t)$ for the process $Y$ of the
lemma. We need\vspace*{-2pt} to show that $\int g \,d\nu_t \xrightarrow{t\to\infty}
\int g \,d\mu$ for any fixed bounded and continuous function
$g\dvtx[-1,1]^\Z\to\R$.

If we start the $N$-type TASEP with an i.i.d. product measure initial
distribution then its distribution converges to an ergodic stationary
measure with the same one-dimensional marginal. (This can be shown by
standard coupling arguments introduced by Liggett; see, e.g.,
\cite{Liggettcouple} or~\cite{Liggett}, Chapter 8.)

Using Lemma~\ref{L:projection} and Corollary~\ref{C:characterization} it
follows that for any nondecreasing step function $F$ on $[-1,1]$ the
process $\{F(Y_n(t))\}_n$ converges in distribution to $\{F(U_n)\}_n$.

For an integer $M$ let $F_M(x) = \frac{\lfloor Mx \rfloor}{M}$, which
maps $[-1,1]$ to $\{i/M, i\in[-M$, $M-1]\}$. Define the operator
$F^\otimes_M$ on configurations, as the operator that applies $F_M$ to
each coordinate: $F^\otimes_M(\eta)_n = F_M(\eta_n)$. Since $g$ is
continuous we can select $M$ such that $\|g - g\circ F^{\otimes}_M
\|_\infty\leq\eps$. By the triangle inequality we have
\[
\biggl|\int g \,d\nu_t - \int g \,d\mu\biggr| \leq2\eps+
\biggl|\int g \circ F^{\otimes}_M \,d\nu_t - \int g\circ F^{\otimes}_M
\,d\mu\biggr|,
\]
and $g\circ F^\otimes_M$ is $g$ applied to a TASEP with finitely many types,
so it can be made smaller than $\eps$ by taking $t$ large enough.
\end{pf}

\section{Two-dimensional marginals of the TASEP speed process}
\label{S:joint}

The key tool for analyzing the joint densities of the speed process is
Corollary~\ref{C:projection_to_N}. This states that if the speed
process is
monotonously projected into $\{1,\ldots,k,\break k+1\}$, then the result is the
stationary measure of the multi-type TASEP with suitable densities. In the
TASEP, the latter is given in terms of the multi-line process (see
Section~\ref{S:stat_tasep}). More explicitly, we will use the following
projections, to which we refer as canonical projections. Let
$\bx=(x_1,\ldots,x_k)$ be an increasing sequence taking values in $[0,1]$,
with the conventions that $x_0=0$ and $x_{k+1}=1$. Define $F\dvtx[-1,1]
\to
\{1,\ldots,k,k+1\}$ by
\[
F(u) = F_{\bx}(u) = \min\{i\dvtx\hat{u} < x_i\}\qquad \mbox{where }\hat
u =
\frac{1+u}{2}.
\]
Note that if $u$ is uniform on $[-1,1]$, then $F(u)=i$ with probability
$x_i-x_{i-1}$. Let $V_i = F(U_i)$, so each $V_i$ has distribution
controlled by the $x$'s. It is not hard to see that the $\sigma$-field
generated by $V_1,\ldots,V_k$ (or any $k$ fixed indices) for all possible
$\bx$'s is the same as the $\sigma$-field of $U_1,\ldots,U_k$.

The scheme of our argument should now be clear. The distribution of $V$ is
given by a multi-line process, and can be computed explicitly. Considering
the resulting probabilities as functions of $\bx$ allows us to recover the
joint density of the corresponding speeds. This last step is done by taking
suitable derivatives w.r.t. $x_i$'s to get the density. In order to find
the joint density of $k$ particles we work with the $k$-line process. In
this section we use this approach to prove results about two-dimensional
marginals of $\mu$. We prove Theorem~\ref{T:joint2} which gives the joint
distribution of $(U_0, U_1)$ and generalize this result for the joint
distribution of any two speeds. In the next section we give some results
for higher-dimensional marginals.

\subsection{Two consecutive speeds: $U_0,U_1$}\label{subs:U0U2}

\mbox{}

\begin{pf*}{Proof of Theorem~\ref{T:joint2}}
We compute the probability that $V_1=2$ and $V_0$ is each of $1,2,3$
(recall that as the highest class particles, $3$'s are equivalent to
holes). The queue of the two line process is a single, simple queue, so
indices are not needed. In order to have a second class particle at
position $1$ we need an unused service. This means the queue must be
empty: $Q(2)=0$, and there must be a particle at the bottom line but not
at the top line in position $1$. The intersection of these events has
probability $x_2-x_1$ (as this is the density of second class particles).
More importantly, they depend only on the two-line configuration in
positions $\{1,\ldots,\infty\}$. Since on this event the queue is also
empty at position 1, the class $V_0$ depends only on the two-line
configuration at position 0.

In particular, to get a first class particle, $V_0=1$, the only
possibility is to also have particles in both lines in position 0. This
leads to
\[
\P(V_0=1,V_1=2 ) = \P\biggl( \config{2}{\D&\C\\\D&\D}^\varnothing
\biggr) = x_1 x_2 (x_2-x_1).
\]
We shall also denote this probability by $\mu_\bx(1,2)$ for compactness,
as this is the probability of seeing consecutive particles of classes
$1,2$ in the stationary measure~$\mu_\bx$. Similarly we have
\begin{eqnarray*}
\mu_\bx(2,2 ) &=& \P\biggl(\config{2}{\C&\C\\\D&\D}^\varnothing
\biggr) = x_1 \bar{x}_2 (x_2-x_1), \\
\mu_\bx(3,2 ) &=& \P\biggl(\config{2}{\B&\C\\\C&\D}^\varnothing
\biggr) = \bar{x}_1 (x_2-x_1).
\end{eqnarray*}
Here, $\B$ indicates no restriction on the top line in that position and
$\bar y=1-y$.

To calculate the densities of the two speeds we find, for example,
\[
\P(U_0 < 2x_1-1 < U_1 < 2x_2-1) = \mu_\bx(1,2) = x_1 x_2 (x_2-x_1).
\]
Thus to find the density at $(u_0,u_1)$ for $u_0<u_1$ we need to take
derivatives w.r.t. $x_2$ and $x_1$, and set $x_2=(1+u_1)/2$,
$x_1=(1+u_0)/2$. Remembering the Jacobians ($1/2$) we find
\begin{eqnarray}
\P(U_0\in d u_0, U_1\in d u_1) =
\biggl( \frac12 \,\partial_{x_1} \biggr)
\biggl( \frac12 \,\partial_{x_2} \biggr) \mu_\bx(1,2)
= \frac{u_1-u_0}{4} \,d u_0 \,d u_1
\nonumber\\
&&\eqntext{\mbox{for $u_0<u_1$.}}
\end{eqnarray}

Similarly, to find the density at $(u_0,u_1)$ for $u_0>u_1$ noting that
the Jacobians now have reversed signs we find
\begin{eqnarray}
\P(U_0\in du_0, U_1\in du_1) =
\biggl( -\frac12 \,\partial_{x_1} \biggr)
\biggl( -\frac12 \,\partial_{x_2} \biggr) \mu_\bx(3,2)
= \frac14 \,d u_0 \,d u_1\nonumber\\
&&\eqntext{\mbox{for $u_0>u_1$.}}
\end{eqnarray}

Finally, to find the (singular) density along the diagonal, consider
$\mu_\bx(2,2)$ and let $x_2,x_1\to\frac{1+u}{2}$. We have
\[
\P(U_0,U_1\in du)
= \frac12 \lim_{x_1,x_2\to(1+u)/2} \frac{\mu_\bx(2,2)}{x_2-x_1}
= \frac{1-u^2}{8} \,du.
\]
\upqed
\end{pf*}

\subsection{Two distant speeds: $U_0,U_k$}\label{subs:U0Uk}

The two line process also yields formulae for the joint density of two
distant particles. However, the result is not as compact as for the
case of
two consecutive particles.
\begin{theorem}\label{T:2distant}
For any $k>0$ we have:
\begin{itemize}
\item The joint density of $U_0,U_k$ on $\{U_0>U_k\}$ is $1/4$ [so
$\P(U_0>U_k)=1/2$].
\item On $\{U_0<U_k\}$ the density is a polynomial 
of degree
$2k-1$.
\item On the diagonal $\{U_0=U_k\}$ the density is a polynomial
of degree $2k$. As $k\to\infty$, the density on the diagonal
$\{(u,u)\dvtx
|u|\le1\}$ is asymptotically $\sqrt{\frac{1-u^2}{16\pi k}}$.
\end{itemize}
\end{theorem}

It is possible to prove exponential convergence of the density on
$\{U_0<U_k\}$ to $1/4$, though we do not pursue that direction here.
The fact
that as $k\to\infty$ the distributions of $U_0$ and $U_k$ become
independent follows from ergodicity, or can be read from \eqref{eq:dist_up}
below.

The theorem follows easily from the next two lemmas. Let $\{S_n\}$ be a
random walk with steps in $\{1,-1, 0\}$ with probabilities
$\{p_+,p_-,p_0\}$, and consider the maximum process $M_n = \max_{i\le
n} S_i$.
\begin{lemma}\label{lemma6.2}
Fix $0\le x<y\le1$, and let $S_n,M_n$ be as above with
\[
p_+ = x\bar{y},\qquad p_- = \bar{x}y,\qquad p_0 = x y + \bar{x} \, \bar{y}.
\]
Then we have the following:
%
%
\begin{eqnarray}
\label{eq:dist_down}
\P(x < \hat U_k < y < \hat U_0) &=& (y-x)\bar{y},
\\
\label{eq:dist_eq}
\P(\hat U_0, \hat U_k \in[x,y])
&=& (y-x)\bar{x} y \P(M_{k-1} = 0),
\\
\label{eq:dist_up}
\P(\hat U_0 < x < \hat U_k < y)
&=& (y-x) x y + (y-x)\bar{x} y \P(M_{k-1}>0).
\end{eqnarray}
\end{lemma}

Note that the steps of $S$ are the difference of two Bernoulli random
variables, and therefore $S_j \stackrel{\mathrm{d}}{=}\bin(j,x) - \bin(j,y)$. In particular,
for any fixed $x<y$ we have $S_j\xrightarrow{j\to\infty}\xxrightarrow
{\mathrm{prob}.}-\infty$,
and asymptotically the speeds are independent.
\begin{pf*}{Proof of Lemma~\ref{lemma6.2}}
By Corollary~\ref{C:projection_to_N}, $\P(x < \hat U_k < y < \hat U_0) =
\mu_{x,y}(3,2)$ (where $\mu_{x,y}$ the extremal stationary $3$ type TASEP
with densities $x, y-x,\break 1-y$). Using the\vadjust{\goodbreak} two-line description of
$\mu_{x,y}$ we have $V_0=1, V_k=2$ if and only if we see the two-line
configuration
\[
\config{3}{\B&\cdots&\C\\\C&\cdots&\D}^\varnothing.
\]
Having the hole in the bottom line at position 0 has probability
$\bar{y}$ and this is independent of having a second class particle at
position $k$.

Similarly, to have $\hat U_0,\hat U_k\in[x,y]$ we need the configuration
\[
\config{3}{\C&\cdots&\C\\\D&\cdots&\D}^\varnothing
\]
with intermediate configuration leaving the queue empty at position 1.
Let $S_j$ be the number of particles in the top line in positions
$\{1,\ldots,j\}$ minus the number of particles in the bottom line in those
positions. The condition that the queue ends up nonempty is equivalent
to $\{\max_{0<j<k-1} S_j \ge1\}$. The claim follows.

Finally, the third case follows from the first two since the three
probabilities must add up to $\P(\hat U_k\in[x,y]) = y-x$.
\end{pf*}
\begin{lemma}
Let $S_n,M_n$ be as above with $p_+=p_-$. Then $\P(M_n = 0 )
= \P(S_n \in\{0,-1\})$.
\end{lemma}
\begin{pf}
Reflection at the hitting time of $1$ shows that
\[
\P(M_n>0, S_n\leq0) = \P(M_n>0, S_n\geq2)
= \P(S_n \geq2) = \P(S_n \leq-2).
\]
It follows that
\[
\P(M_n>0) = \P(S_n>0) + \P(M_n>0, S_n\leq0)
= 1 - \P(S_n \in\{0,-1\}).
\]
\upqed
\end{pf}
%
%
\begin{pf*}{Proof of Theorem~\ref{T:2distant}}
The case $U_0>U_k$ is just the double derivative of \eqref{eq:dist_down}.

For the case $U_0=U_k$, note from \eqref{eq:dist_eq} that the density
along the diagonal~is
\[
\lim_{x,y\to\hat{u}} \frac{\P(\hat U_0,\hat U_k \in[x,y])}{2(y-x)}
= \frac{1-u^2}{8} \P(M_{k-1}=0),
\]
where $M_{k-1}$ is the maximum of a symmetric random walk with
$p_+=p_-=x\bar{x}$. Using the prior lemma, since $p_+=p_-$ we get
\[
\P(M_{k-1}=0) = \P(S_{k-1} \in\{0,-1\}).
\]
This is clearly polynomial. Using the local central limit theorem,
$\P(S_{k-1}=a) \sim\frac{1}{\sqrt{4\pi x\bar{x} k}}$ for any
$a\in\{0,-1\}$, and our claims follow.

For the case $U_0<U_k$, taking derivatives of \eqref{eq:dist_up} shows
that the density is polynomial as claimed.
\end{pf*}

\section{Multiple speeds} \label{S:mult}

In this section we will prove some results about the joint distribution of
more than two speeds. In principle, any finite-dimensional marginal of the
distribution can be derived from Theorem~\ref{T:stationary} along the\vadjust{\goodbreak}
same lines
as used above for the joint distribution of $U_0,U_1$. This gives the joint
distribution in terms of the stationary measure of the multiple queue
system. Some aspects of the joint distribution have particularly nice
formulae, and we proceed to present some of these:
\begin{longlist}[(2)]
\item[(1)] The next subsection determines the probability that out of the first
$n$ particles a given one is the fastest.
\item[(2)] The following result shows that the speed of a fast particle is
independent from those of adjacent particles it overtakes. More
precisely, if $c\in[-1,1]$, then conditioned on the event that $U_0 >
c $
and $U_1, \ldots, U_n < c$, the random vector $(U_1,\ldots,U_n)$ and $U_0$
are independent.
\item[(3)] Next, we show that on the event $\{U_0<U_1<\cdots<U_n\}$ there is a
pairwise repulsion between the particles: the density function is given
by $n!$ times a Vandermonde determinant.
\item[(4)] Finally, we give the full description of the joint distribution of
$(U_0,\break U_1,U_2)$. Their distribution is absolutely continuous with respect
to the Lebes\-gue measure on each of the 13 subsets of $[-1,1]^3$
corresponding to a given order of these speeds (these include the cases
where two or three speeds might be equal). In Theorem~\ref{T:3speeds} we
determine the densities on all of these subsets.
\end{longlist}

\subsection{The fastest particle}

As a first example, we compute the probability that particle $i$ will be
the rightmost of $\{1,\ldots,n\}$ for all $t>t_0$. This proves and
generalizes a conjecture of Ferrari, Goncalves and Martin~\cite{FGM} that
the probability of particle $0$ overtaking particles $1$ through $n$ is
$\frac{2}{n+2}$. Note that this is not quite the same as saying that $U_i$
is the maximal of $\{U_1,\ldots,U_n\}$. Due to Lemma~\ref{L:diagswap},
this event
allows $U_i=U_j$ for $j>i$ but not for $j<i$.
\begin{theorem}\label{T:rightmost}
For any $n$ and any $k\in[1,n]$
\[
\lim_{t\to\infty} \P\bigl(X_{k}(t) = \max\{X_{1}(t),\ldots,X_{n}(t) \}
\bigr) = \frac{2n}{(n+k-1)(n+k)}.
\]
\end{theorem}
\begin{lemma}
Let $X \eqd\bin(m,p)$ and $Y \eqd\geom(q)$ be independent binomial and
geometric random variables. Then
\[
\P(Y\le X) = 1 - q (\bar{p} + p q)^m.
\]
\end{lemma}
\begin{pf}
We have
$
\P(Y > X)
= \sum_i {m\choose i} p^i \bar{p}^{m-i} q^{i+1}
= q (\bar{p} + p q )^m
$.
\end{pf}
\begin{pf*}{Proof of Theorem~\ref{T:rightmost}}
Since the index of the rightmost particle (among the set $\{1,\ldots,n\}
$) is
nonincreasing in time, the event in the statement is equivalent to
particle $k$ being the rightmost for all $t>t_0$ for some~$t_0$. By
Lemma~\ref{L:diagswap}, which we prove in Section~\ref{S:ASEP},
particle $i$
eventually passes particle $j$ for $i<j$ if and only if $U_i \geq U_j$.
Thus $k$ will eventually be the rightmost particle of particles
$\{1,\ldots,n\}$ if and only if $U_k > U_i$ for $1\le i < k$ and $U_k
\ge
U_i$ for $k < i \le n$. Call this event $E_k$.\vadjust{\goodbreak}

As an intermediate step we will compute the probability that this happens
and $U_k\in du$ for some $u\in[-1,1]$. Integrating over $u$ will give the
theorem. Fix $\bx=(x_1,x_2)$, $0<x_1<x_2<1$ and consider the event
$E_{k,\bx}$ that for all $i\in[1,n]$ we have that
\[
\hat U_i \in\cases{
[0,x_1], &\quad$i<k$, \cr
[x_1,x_2], &\quad$i=k$, \cr
[0,x_2], &\quad$i>k$.}
\]
Thus $E_{k,\bx}$ says that up to the partition resulting from the vector
$\bx$, the event $E_k$ holds.

Projecting into the $2+1$ type TASEP using $F_{\bx}$, $E_{k,\bx}$ is
mapped to the of event of having $k-1$ first class particles followed by
a second class particle, followed by $n-k$ particles of either class (but
no holes). This requires in positions $1$--$n$ a configuration of the
following form:
\[
\config{7}{\D&\cdots&\D&\C&\B&\cdots&\B\\\D&\cdots&\D&\D&\D&\cdots&\D}^i,
\]
where the first hole in the top line is in position $k$, and the size $i$
of the queue can be no greater than the number of holes in the top line
in positions $\{k+1,\ldots,n\}$. Since the number of holes in the rest of
the top line has the binomial distribution $\bin(n-k,\bar{x}_1)$ and the
queue state is an independent $\geom(
\frac{x_1\bar{x}_2}{\bar{x}_1x_2} )$, we find after simplifying
that
\[\hspace*{-5pt}
\P(E_{k,\bx}) = x_1^{k-1} \bar{x}_1 x_2^n \P\biggl(
\geom\biggl( \frac{x_1\bar{x}_2}{\bar{x}_1x_2} \biggr) \le
\bin(n-k,\bar{x}_1) \biggr)
= x_1^{k-1} \bar{x}_1 x_2^n - x_1^n x_2^{k-1} \bar{x}_2
\]
(noting that $\bar{q} + p q$ of the previous lemma simplifies to
$x_1/x_2$).

Taking a limit as $x_2,x_1\to y$ we find
\[
\P(E_k, \hat U_k \in d y) = \lim_{x_2,x_1\to y}
\frac{\P(E_{k,\bx})}{x_2-x_1}
= y^{n+k-2}\bigl((n+1-k) - (n-k)y \bigr) \,dy.
\]
Finally, integrating over $y\in[0,1]$ gives
\[
\P(E_k) =\int_0^1 y^{n+k-2}\bigl((n+1-k) - (n-k)y \bigr) \,dy
= \frac{2(n+1)}{(n+k-1)(n+k)}.\quad
\]
\upqed
\end{pf*}

\subsection{Independence when swapped}

The following result shows that the speed of a fast particle is independent
of speeds of adjacent particles that it overtakes.
\begin{lemma}\label{L:adding_fast_one}
Fix $c\in[-1,1]$ and a measurable set $A\subset[-1,c]^n$. Then we have
\[
\mu\bigl( U_0 > c | (U_1,\ldots,U_n)\in A\bigr) = \mu(U_0>c).
\]
Furthermore, conditioned on $U_0>c$ and $(U_1,\ldots,U_n)\in A$ we have
that $U_0$ is uniform on $[c,1]$.
\end{lemma}
\begin{pf}
Since products of intervals span the $\sigma$-field, it suffices to prove
the analogous statement for the $M$-type TASEP (in fact $M=n+1$ is
enough). Consider a TASEP measure $\mu_\bx$ where holes have density
$1-\hat{c}$, so that speeds greater than $c$ correspond to holes. We need
to show that for any classes $i_1,\ldots,i_n < M$
%
%
\begin{equation}\label{eq:cond_ind}
\mu_\bx(V_0=M | V_1=i_1,\ldots,V_n=i_n) = \mu_\bx(V_0=M).
\end{equation}
To show this we consider the multi-line process. There the classes of
$V_1,\ldots,V_n$ are determined by the lines in positions $[1,\infty)$. On
the other hand, $V_0=M$ requires only that $B_M(0)=0$, hence the
independence.

To get the second claim, note that $\mu(U_0>c)=\frac{1-c}{2}$ and that
\eqref{eq:cond_ind} also applies (with the same set $A$) for any $c'>c$.
\end{pf}
\begin{coro}\label{C:reversed_order}
The $(U_1,\ldots,U_n)$-marginal of $\mu$ has a constant density function
$2^{-n}$ on the set $\{U_1 > \cdots> U_n\}$.
\end{coro}
\begin{pf}
The events that the speeds are in small intervals around the $u_i$'s are
independent.
\end{pf}

\subsection{Repulsion when unswapped}

Here we derive the density function of the $(n+1)$-dimensional marginal of
$\mu$ on the event $\{U_0<\cdots<U_n\}$. The result is given in terms of
a Vandermonde determinant defined by
\[
\Delta_{a,b}(\bx) = \prod_{a\leq i < j \leq b} (x_j-x_i).
\]
We start with a simple lemma about these determinants.
\begin{lemma}\label{L:vandermonde}
Let $x_0 < \cdots< x_n$. Then
\[
\Delta_{0,n}({\bx}) = n! \int_{x_{i-1}<y_i<x_i} \Delta_{1,n}(y)
\prod_{i=1}^n d y_i.
\]
\end{lemma}
\begin{pf}
We use the standard fact that $\Delta(y)$ is the determinant of the
Vandermonde matrix: $\Delta_{1,n}(y) = \det
(y_i^{j-1})_{i,j=1}^n$. Since the determinant is is linear in
the rows and each $y_i$ appears in a single row, we can integrate row by
row to find
\begin{eqnarray*}
\int_{x_{i-1}<y_i<x_i} \Delta_{1,n}(y) \prod_{i=1}^n d y_i
&=& \det\int_{x_{i-1}<y_i<x_i}
(y_i^{j-1})_{i,j=1}^n \prod_{i=1}^n d y_i \\
&=& \det\biggl( \frac{x_{i-1}^j-x_i^j}{j} \biggr)_{i,j=1}^n
= \frac{1}{n!} \det M,
\end{eqnarray*}
where $M=( x_{i-1}^j-x_i^j )_{i,j=1}^n$. Extend $M$ to an
$(n+1)\times(n+1)$ matrix $M'$ by
\[
M' = \left(
\begin{array}{c|c}
1 & \matrix{ x_0 & \cdots& x_0^n } \cr
\hline
\matrix{ 0 \cr\vdots\cr0 } & M
\end{array}
\right).
\]
Clearly $\det M = \det M'$. However, by sequentially adding each row to
the one below it we find $\det M' = \det(x_i^{j-1})_{i,j=0}^n
= \Delta_{0,n}({\bx})$, completing the proof.
\end{pf}
\begin{lemma}\label{L:preserved_order}
Let $0 = x_0 < x_1 < \cdots< x_n < x_{n+1} = 1$, and $\mu_{\bx}$ be the
corresponding $n+1$ type TASEP stationary measure. Let $Q_n$ be the
probability that all queues are empty at any specific location of the $n$
line process. We have the following:
\begin{longlist}[(2)]
\item[(1)]
$\mu_{\bx}(2,\dots,n) = \mu(\hat{U_i}\in[x_{i-1},x_i]$ for all $i\in[2,n]) =
\Delta_{1,n}({\bx})$,
\item[(2)] $\mu_{\bx}(1,\dots,n) = \mu(\hat{U_i}\in[x_{i-1},x_i]$ for all $i \in[1,n]) =
\Delta_{0,n}({\bx})$,
\item[(3)] The density of $\hat U_1,\ldots,\hat U_n$ on the event $U_1 <
\cdots
< U_n$ is $n! \Delta_{1,n}(\hat u)$;
\item[(4)] $Q_n = \frac{\Delta_{1,n}({\bx})} {\prod_{i=1}^n x_i^{i-1}
\bar{x}_i^{n-i} }$.
\end{longlist}
\end{lemma}
\begin{pf}
The proof is by induction on $n$. For $n=1$, claims (1) and (4) are trivially
true, and (2), (3) hold since the speeds are uniformly distributed.

The key observation is that the only $n$-line configuration giving
particles of classes $1,\ldots,n$ is
\[
\config{6}{
\D&\C&\C&\cdots&\C&\C\\
\D&\D&\C&\cdots&\C&\C\\
\D&\D&\D&\cdots&\C&\C\\
\vdots&\vdots&\vdots&\ddots&\vdots&\vdots\\
\D&\D&\D&\cdots&\D&\C\\
\D&\D&\D&\cdots&\D&\D}^{\varnothing,\ldots,\varnothing}
\]
(with all queues empty). Since the queue state is independent of the
configuration in these positions, we find
\[
\mu_{\bx}(1,\ldots,n) = Q_n \prod_{i=1}^n x_i^{i} \bar{x}_i^{\hspace*{0.1pt}n-i}.
\]
This implies equivalence of claims (2) and (4).

Similarly, the only configuration giving particles of types $2,\ldots,n$
is
\[
\config{5}{
\C&\C&\cdots&\C&\C\\
\D&\C&\cdots&\C&\C\\
\D&\D&\cdots&\C&\C\\
\vdots&\vdots&\ddots&\vdots&\vdots\\
\D&\D&\cdots&\D&\C\\
\D&\D&\cdots&\D&\D}^{\varnothing,\ldots,\varnothing}.
\]
This implies equivalence of claims (1) and (4) [since
$\Delta_{0,n}({\bx})=\Delta_{1,n}({\bx})\prod x_i$].\vadjust{\goodbreak}

Next, we argue that claims (2) and (3) are equivalent. Claim (2) follows
from (3)
by Lemma~\ref{L:vandermonde}. Claim (2) also implies claim (3), since the
density is the multiple derivative $\prod_{i=1}^n \frac\partial
{\partial
x_i}$ of the probability of claim (2).

Thus for any given $n$, the four claims are all equivalent. To complete
the proof (by induction) we note that claim (3) for a given $n$ implies
claim (1) for $n+1$. This also follows from Lemma~\ref{L:vandermonde} in the
same way as claim (2).
\end{pf}

\subsection{Joint densities for 3 consecutive particles}

This section contains the complete description of the joint
distribution of
$(U_0,U_1,U_2)$. The distribution is absolutely continuous with respect to
the Lebesgue measure on each of the 13 subsets of $[-1,1]^3$ corresponding
to a given order of these speeds (these include the cases where two or all
three speeds might be equal). In Theorem~\ref{T:3speeds} we determine the
densities on all of these subsets.
\begin{theorem}\label{T:3speeds}
The joint distribution of $U_0,U_1,U_2$ is given by Table~\ref{table2},
%
\begin{table}
\caption{Joint densities of $(U_0, U_1, U_2)$ according to their relative order}\label{table2}
\begin{tabular*}{240pt}{@{\extracolsep{\fill}}lc@{}}
\hline
\textbf{Order} & \textbf{Density} \\
\hline
$u_0 < u_1 < u_2$ & $\frac{3}{32} (u_2-u_1) (u_1-u_0) (u_2-u_0)$ \\[2pt]
$u_0 < u_2 < u_1$ & $\frac{1}{32} (u_2-u_0) (2+4u_1-3u_2-3u_0)$ \\[2pt]
$u_1 < u_0 < u_2$ & $\frac{1}{32} (u_2-u_0) (2+3u_2+3u_0-4u_1)$ \\[2pt]
$u_1 < u_2 < u_0$ & $\frac{1}{8} (u_2-u_1)$ \\[2pt]
$u_2 < u_0 < u_1$ & $\frac{1}{8} (u_1-u_0)$ \\[2pt]
$u_2 < u_1 < u_0$ & $\frac{1}{8}$ \\[2pt]
$u_0 = u_1 < u_2$ & $\frac{1}{64} (u_2-u_1)(1-u_1^2)(2+3u_2-u_1)$ \\[2pt]
$u_0 < u_1 = u_2$ & $\frac{1}{64} (u_1-u_0)(1-u_1^2)(2-3u_0+u_1)$ \\[2pt]
$u_1 < u_0 = u_2$ & $\frac{1}{16} (u_2-u_1)(1-u_2^2)$ \\[2pt]
$u_0 = u_2 < u_1$ & $\frac{1}{16} (u_1-u_0)(1-u_0^2)$ \\[2pt]
$u_1 = u_2 < u_0$ & $\frac{1}{16} (1-u_1^2)$ \\[2pt]
$u_2 < u_0 = u_1$ & $\frac{1}{16} (1-u_1^2)$ \\[2pt]
$u_0=u_1=u_2$ & $\frac{1}{32} (1-u_0^2)^2$ \\
\hline
\end{tabular*}
\end{table}
arranged according to their relative order.
\end{theorem}
\begin{pf}
Fix $0<x_1<x_2<x_3<1$. Define $F=F_{\bx}$ as above, and $V_i=F(U_i)$. To
calculate the densities of the various simplices and facets, we calculate
partly the distribution of $V$, and take suitable derivatives and limits.
It is interesting to note that there are several possible class
configurations for each case. For example, the case $\{U_0<U_1<U_2\}$ can
be deduced from each of $\mu_{\bx}(1,2,3)$, $\mu_{\bx}(1,2,4)$,
$\mu_{\bx}(1,3,4)$ and $\mu_{\bx}(2,3,4)$. Careful choice of the cases
to consider can simplify the computations significantly.\vadjust{\goodbreak}

Not all cases need to be worked out. Space-class symmetry reduces several
cases to others. Theorem~\ref{T:joint2}, Lemmas \ref
{L:adding_fast_one} and
\ref{L:preserved_order} and Corollary~\ref{C:reversed_order} imply
several cases. Thus even though all 13 cases can be computed using this
method, only 4 are essentially new and proved below.

Table~\ref{table3} summarizes the proofs for the 13 weak orders of
$U_0,U_1,U_2$. Here, $\Delta({\bx})=\Delta_{1,3}({\bx})$.

%
%
\begin{table}
\caption{Main ingredients for the proofs of the joint densities in Theorem \protect\ref{T:3speeds}}\label{table3}
\begin{tabular*}{\tablewidth}{@{\extracolsep{\fill}}lccl@{}}
\hline
\textbf{Order} & $\bolds{V}$ & $\bolds{\mu_{\bx}(V)}$ & \multicolumn{1}{c@{}}{\textbf{Remarks}} \\
\hline
$U_0 < U_1 < U_2$ & $1,2,3$ & $x_1 x_2 x_3 \Delta({\bx})$ &
\mbox{Lemma~\ref{L:preserved_order}} \\
[2pt]
$U_0 < U_2 < U_1$ & $2,4,3$ &
$(\bar{x}_1+\bar{x}_2) \bar{x}_3 \Delta({\bx})$ & New \\
[2pt]
$U_1 < U_0 < U_2$ & $2,1,3$ & $x_1(x_2+x_3) \Delta({\bx})$ &
Space-class symmetry \\
[2pt]
$U_1 < U_2 < U_0$ & $4,2,3$ & $\bar{x}_3 \Delta({\bx})$ &
Theorem~\ref{T:joint2}, Lemma~\ref{L:adding_fast_one}\\
[2pt]
$U_2 < U_0 < U_1$ & $2,3,1$ & $x_1 \Delta({\bx})$ &
Space-class symmetry \\
[2pt]
$U_2 < U_1 < U_0$ & $3,2,1$ & $x_1(x_2-x_1)(x_3-x_2)$ &
Corollary~\ref{C:reversed_order} \\
[6pt]
$U_0 = U_1 < U_2$ & $2,2,3$ & $\bar{x}_1 x_2 x_3 \Delta({\bx})$ & New
\\
[2pt]
$U_0 < U_1 = U_2$ & $2,3,3$ & $\bar{x}_1 \bar{x}_2 x_3 \Delta({\bx})$ &
Space-class symmetry \\
[2pt]
$U_1 < U_0 = U_2$ & $3,2,3$ & $\bar{x}_2 x_3 \Delta({\bx})$ & New \\
[2pt]
$U_0 = U_2 < U_1$ & $2,3,2$ & $\bar{x}_1 x_2 \Delta({\bx})$ &
Space-class symmetry \\
[2pt]
$U_1 = U_2 < U_0$ & $4,2,2$ & $\bar{x}_1 x_2 \bar{x}_3 (x_2-x_1)$ &
Theorem~\ref{T:joint2}, Lemma~\ref{L:adding_fast_one} \\
[2pt]
$U_2 < U_0 = U_1$ & $3,3,1$ & $x_1 \bar{x}_2 x_3 (x_3-x_2)$ &
Space-class symmetry \\
[6pt]
$U_0 = U_1 = U_2$ & $2,2,2$ & $\bar{x}_1^2 x_2^2 (x_2-x_1)$ &
New; Theorem~\ref{T:convoys} \\
\hline
\end{tabular*}
\end{table}

The case $\{U_0<U_1<U_2\}$ is a special case of Lemma~\ref{L:preserved_order},
while the case $\{U_2<U_1 < U_0\}$ is a special case of
Corollary~\ref{C:reversed_order}. The cases $\{U_1 < U_2 < U_0\}$ and
$\{U_1 = U_2
< U_0\}$ follow from joint distribution of $U_1,U_2$ (Theorem~\ref{T:joint2})
together with Lemma~\ref{L:adding_fast_one}. Each of the five
cases $\{U_2<U_0<U_1\}$, $\{U_1<U_0<U_2\}$, $\{U_0<U_1=U_2\}$,
$\{U_0=U_2<U_1\}$ and $\{U_2<U_0=U_1\}$ follows by space-class symmetry
(Proposition~\ref{P:spacesym}) from the cases $\{U_1<U_2<U_0\}$, $\{
U_0<U_2<U_1\}$,
$\{U_0=U_1<U_2\}$, $\{U_1<U_0=U_2\}$ and $\{U_1=U_2<U_0\}$, respectively.

It therefore remains to prove just 4 cases: $\{U_0<U_2<U_1\}$,
$\{U_0=U_1<U_2\}$, $\{U_1<U_0=U_2\}$ and $\{U_0=U_1=U_2\}$.

For the case $\{U_0 < U_2 < U_1\}$, we compute $\mu_{\bx}(2,4,3)$.
The only 3 line configurations that give these types are
\[
\config{3}{\B&\C&\C\\\C&\D&\C\\\D&\C&\D}^{\varnothing,\varnothing}
\quad\mbox{and}\quad
\config{3}{\C&\C&\C\\\D&\B&\C\\\D&\C&\D}^{\varnothing,\varnothing}.
\]
Therefore
\begin{eqnarray*}
\mu_{\bx}(2,4,3)
&=& \bar{x}_1^2 x_2\bar{x}_2 x_3^2 \bar{x}_3
(\bar{x}_1+\bar{x}_2) \mu_{\bx}(\mbox{empty queues}) \\
&=& \bar{x}_1^2 x_2\bar{x}_2 x_3^2 \bar{x}_3
(\bar{x}_1+\bar{x}_2) \frac{\Delta_{1,3}({\bx})}{\bar{x}_1^2
x_2\bar{x}_2 x_3^2} \\
&=& \bar{x}_3 (\bar{x}_1+\bar{x}_2) {\Delta_{1,3}({\bx})}.
\end{eqnarray*}
Taking derivatives we find the density of $\hat U_0,\hat U_1,\hat U_2$ in
the domain $\{U_0<U_2<U_1\}$ is
\[
\frac{-\partial}{\partial x_3}\bigg|_{x_3=\hat u_1}
\frac{-\partial}{\partial x_2}\bigg|_{x_2=\hat u_2}
\frac{-\partial}{\partial x_1}\bigg|_{x_1=\hat u_0}
\mu_{\bx}(2,4,3) =
(\hat u_1 - \hat u_2) (2 + 4 \hat u_1 - 3 \hat u_0 -3 \hat u_2).
\]
A linear change of variables gives the formula in terms of $u_1,u_2,u_3$.

For the case $\{U_0=U_1<U_2\}$, we consider $\mu_{\bx}(2,2,3)$. The only
three-line configuration giving this result is
\[
\config{3}{\C&\C&\C\\\D&\D&\C\\\D&\D&\D}^{\varnothing,\varnothing}.
\]
Thus
\[
\mu_{\bx}(2,2,3)
= x_3^3 x_2^2\bar{x}_2\bar{x}_1^3 \mu_{\bx}(\mbox{empty queues})
= \bar{x}_1 x_2 x_3 \Delta_{1,3}({\bx}).
\]
Taking a derivative w.r.t. $x_3$ and letting $x_2\to x_1$ gives the
density of the $\hat U_i$'s to be
\[
\lim_{x_1,x_2\to\hat u_0} \frac{1}{x_2-x_1}
\frac{-\partial}{\partial x_3}\bigg|_{x_3=\hat u_2}
\mu_{\bx}(2,2,3)
= \hat u_0 \bar{\hat u}_0 (\hat u_2 - \hat u_0)(3 \hat u_2 - \hat u_0).
\]
As above, a change of variables gives the claim.

For the case $\{U_1<U_0=U_2\}$ we consider $\mu_{\bx}(3,2,3)$. The three-line
configurations giving these classes are of the form
\[
\config{3}{\B&\C&\C\\\C&\D&\C\\\D&\D&\D}^{\varnothing,\varnothing}
\]
and therefore
\[
\mu_{\bx}(3,2,3) = \bar{x}_2 x_3 \Delta_{1,3}({\bx}).
\]

Finally, the case $\{U_0=U_1=U_2\}$ is related to the convoys studied in
Section~\ref{S:convoys}. Indeed, the formula follows from the density of
$U_0,U_1$ and the result that convoys are renewal processes. A more
direct approach follows. As there are no third-class particles in this
case, we will use the projection into the $2+1$ type TASEP using only
$x_1,x_2$ (or equivalently, $x_3=x_2$). The only two-line configuration
giving classes $(2,2,2)$ is
\[
\config{3}{\C&\C&\C\\\D&\D&\D}^{\varnothing}
\]
and therefore
\[
\mu_{\bx}(2,2,2) = \bar{x}_1^2 x_2^2 (x_2-x_1).
\]
Dividing by $x_2-x_1$ and taking a limit $x_2\to x_1$ gives the density
$x_1^2\bar{x}_1^2$.
\end{pf}

\section{Convoys} \label{S:convoys}

The convoy phenomenon is the fact that even though each particle's
speed is
uniform on $[-1,1]$, any two particles have positive probability of having
equal speeds. Indeed, a.s. there will be infinitely many particles with
the same speed as any given particle. We refer to such sets of
particles as
convoys. Thus $\Z$ is partitioned in some translation invariant way into
disjoint infinite convoys.

Let $C_k = \{n \dvtx U_n=U_k\}$ denote the convoy of particle $k$, that
is, all
particles with the same speed as $k$. We will restrict ourselves here to
the study of a single convoy, though the multi-line description of the
multi-type stationary distribution can in principle be used to understand
the joint distribution of several convoys.
\begin{pf*}{Proof of Theorem~\ref{T:convoys}}
Partition the particles into three classes, with thresholds
$\bx=(u,u+\eps)$. The stationary measure $\mu_\bx$ has particles of
classes $1,2,3$ with respective densities $u,\eps,1-u-\eps$. It is known
that the second class particles form a renewal process. The key to the
proof is (as above) to condition on $\hat U_0\in[u,u+\eps]$ and let
$\eps\to0$.

Consider the two line process giving $\mu_\bx$, and let $T_k,S_k$ be the
counting functions of particles in the top and bottom lines, respectively,
so that $T_k$ is the number of particles in $(0,k]$ in the top line. We
may extend $S,T$ to negative $k$ by having $S_k$ be minus the number of
particles in $(-k,0]$ and similarly for $T_k$. It is clear that
$\{S_k\},\{T_k\}$ are random walks with $\{0,1\}$ steps with
$\P(S_{k+1}-S_k=1) = u+\eps$ and $\P(T_{k+1}-T_k=1) = u$. Let
$V\in\{1,2,3\}^\Z$ denote the resulting configuration with the stationary
distribution with these densities.

The two-line collapsing procedure implies the identity
\[
\{V_1=2\} = \Bigl\{S_1=1, T_1=0, \min_{k>0} S_k-T_k > 0\Bigr\}
= \Bigl\{\min_{k>0} S_k-T_k > 0\Bigr\}
\]
(since $S_0=T_0=0$). Further, $V_k=2$ if and only if $S_k-S_{k-1}=1,
T_k-T_{k-1}=0$ and $\min_{\ell\geq k} S_\ell-T_\ell= S_k-T_k$. This
suggests looking at the random walk $R_k=S_k-T_k$, with steps with
distribution
\[
\P(R_{k+1}-R_k = x) = \cases{
\bar{u} (u+\eps), &\quad$x=1$, \cr
u(u+\eps) + \bar{u} \,\bar{u+\eps}, &\quad$x=0$, \cr
u\bar{u+\eps}, &\quad$x=-1$.}
\]

Having the second class particle at $1$ implies that $R$ stays positive,
while its drift is $O(\eps)$. As $\eps\to0$ the distribution of $R$
converges (in the product topology for sequences) to a random walk
conditioned to stay positive for all $n>0$ with step distribution
\[
\P(R_{k+1}-R_k = x) = \cases{
u \bar{u}, &\quad$x=\pm1$, \cr
u^2+\bar{u}^2, &\quad$x=0$.}
\]
Thus $R$ is a lazy simple random walk, and the only effect of $u$ is
through the probability of making a nonzero move.\vadjust{\goodbreak} Having a second class
particle at $1$ does not depend on values of $R_n$ for $n<0$, and this is
also the case in the limit as $\eps\to0$.

This random walk conditioned to stay positive will a.s. tend to
$\infty$
as $n\to\infty$. Furthermore, if we take $u=\hat U_0$ then as $\eps\to0$
the second class particles are exactly at $k$ with $U_k=U_0$. In
particular, the convoy $C_1$ is equal in law to the times of the last
visits of $R$ to any value
\[
C_1 = \{ n \dvtx m\geq n\implies R_m\geq R_n \}.
\]

The claim that the convoys are renewal processes follows either from the
corresponding fact about the times of last visits of $R$ conditioned to
remain positive, or from the fact that for any $\eps>0$ the second class
particles form a renewal process.

If the random walk were just a simple random walk (not lazy) then the
probability of having a jump of length $2k+1$ (as even lengths are
impossible) would be $p_{2k+1} = 2^{-(2k+1)} \frac{1}{k+1}{2k\choose k}$.
The laziness of the random walk implies that the distance from a particle
to the next in a convoy with speed $u$ is a sum of $K$ geometric random
variables with mean $1/(2u\bar{u})$ where $\P(K=2k+1)$ is as above. In
particular, $\P(\mathrm{dist}=m) \asymp\frac{c}{u\bar{u}
m^{3/2}}$.\vspace*{-2pt}
\end{pf*}
\begin{example}
Consider $\P(U_0=U_1=\cdots=U_n)$. The probability that all these speeds
are in some infinitesimal $du$ is
\[
\P(\hat U_0,\ldots,\hat U_n \in du) = (u\bar{u})^n \,du.
\]
(This can be seen easily from the corresponding density $u\bar{u}\,du$ for
two particles and the renewal property.) Integrating gives
\[
\P(U_0=\cdots=U_n) = \frac{n!^2}{(2n+1)!}.\vspace*{-2pt}
\]
\end{example}

\section{Joint distribution---ASEP} \label{S:ASEP}

We present two variations of our argument. The first is restricted to
considering the probability that two adjacent particles are unswapped at
large time. This event is roughly equivalent to $\{U_0<U_1\}$, with some
contribution from $\{U_0=U_1\}$.

The second variation came from an attempt to extract the complete joint
distribution of two speed. For the ASEP it is less successful than form the
TASEP, and is also conditional on a.s. existence of the speeds process.\vspace*{-2pt}

\subsection{Swap probabilities}

The key to our analysis of swap probabilities in the ASEP is to double
count swaps happening until time $t$. Let $R(t)$ be the expected number of
particles $j>0$ that are swapped with 0 at time $t$,
that~is,\looseness=-1
\[
R(t) = \E\# \{j>0 \dvtx X_{0}(t) > X_{j}(t)\}.
\]\looseness=0

Recall the time $t$ speed process $U(t)$ is defined by $U_i(t) =
\frac{X_{i}(t)-i}{t}$. Define the empiric time $t$ measure $\nu_t$ by
\[
\nu_t = \frac{1}{t} \sum_i \delta_{{i}/{t}, U_i(t)}.\vadjust{\goodbreak}
\]
The following is equivalent to the standard hydrodynamic limit theorem for
the ASEP started with the Riemann initial condition.
\begin{lemma}
Almost surely $2\rho\nu_t$ converges weakly to the Lebesgue measure on
$\R\times[-\rho,\rho]$.
\end{lemma}

The following simple fact is frequently useful.
\begin{lemma}\label{L:joint_limit}
Let $\mathcal X_1, \mathcal X_2$ be topological spaces and $(X(t),Y(t)),
t\ge0$ be random variables on the product space $\mathcal X_1\times
\mathcal X_2$. Suppose that $X(t)\xrightarrow{t\to\infty}\xxrightarrow
{\mathit{prob}.}x$ and
$Y(t)\xrightarrow{t\to\infty}\xxrightarrow{\mathit{dist}} Y$ where $x\in
\mathcal X_1$ and
$Y$ is an $\mathcal X_2$-valued random variable. Then the joint limit
also holds $(X(t),Y(t)) \xrightarrow{t\to\infty}\xxrightarrow{\mathit{dist}} (x,Y)$.
\end{lemma}

The application in our case involves $X(t)=\nu_t$, which converges in
probability to Lebesgue measure on a stripe (in the space of measures) and
$Y(t)=U_0(t)$ which tends to $U_0$. The conclusion implies that the
hydrodynamic limit also holds conditioned on $U_0$.

The next lemma determines the asymptotic value of $R(t)$.
\begin{lemma}\label{L:R_t_lim}
$R(t) \sim\rho t/3$.
\end{lemma}
\begin{pf}
Particle $0$ has swapped with particle $j>0$ if and only if $X_{j}(t) <
X_{0}(t)$, which can be written as
\[
U_j(t) < U_0(t) - \frac{j}{t}.
\]
It follows that
\[
\frac{R(t)}{t} = \E\bigl[ \nu_t\bigl( \{ (x,y) \dvtx0<x<U_0(t)-y \} \bigr)
\bigr].
\]
Now, Lemma~\ref{L:joint_limit} (see the subsequent discussion) shows
that we
can take a joint limit as $U_0(t)$ converges in distribution to uniform
on $[-\rho,\rho]$, and $\nu_t$ converges weakly in probability to a fixed
measure which is $1/(2\rho)$ times Lebesgue on a strip. Thus
\begin{eqnarray*}
\lim_{t\to\infty}\frac{R(t)}{t}
&=& \E\biggl[ \frac{1}{2\rho} \Leb\bigl( \{ (x,y) \dvtx0<x<U-y, y\ge-\rho
\}
\bigr) \biggr]\\
&=& \E\frac{(U+\rho)^2}{4\rho} = \frac{\rho}{3}.\hspace*{200pt}\qed
\end{eqnarray*}
\noqed\end{pf}

Consider now the following probability (Theorem~\ref{T:symmetry} shows
that the
two definitions are equivalent):
%
%
\begin{equation}\label{eq:Qdef}
Q(t) = \P\bigl(X_{0}(t) < X_{1}(t)\bigr) = \P\bigl(Y_{0}(t) < Y_{1}(t)\bigr).
\end{equation}
$Q(t)$ measures the probability that particles $0$ and $1$ are
unswapped at
time $t$---our present objective.
\begin{lemma}\label{L:Q_mon}
$Q(t)$ is monotone decreasing in $t$.
\end{lemma}
\begin{pf}
Condition on all events except those involving particles $\{0,1\}$, and
denote this $\sigma$-field by $\F_{0,1}$. Recall that $J_{0,1}(t)$
denotes the time 0 and 1 spends next to each other up to time $t$ and
note that $J_{0,1}(t)$ is measurable in $\F_{0,1}$. Then by (\ref{eq:J})
we have $\P(X_{0}(t)<X_{1}(t)|\F_{0,1}) = \bar{p} + p e^{-J_{0,1}(t)}$.
Since $J_{0,1}(t)$ is increasing, $Q(t) = \bar{p} + p\E e^{-J_{0,1}(t)}$
is decreasing.
\end{pf}
\begin{lemma}\label{L:R_t_deriv}
For any $t$ we have $\frac{d}{dt} R(t) = p Q(t) - \bar{p}\bar{Q(t)} =
p +
Q(t) - 1$.
\end{lemma}
\begin{pf}
Let $r_i^+(t)$ [resp., $r_i^-(t)$] be the probability that at time $t$
particle $i$ has a larger indexed particle to its right (resp., left). By
translation invariance these do not depend on $i$. $R(t)$ is the
expectation of a random variable which increases by one with rate $p$ if
the particle at $X_0(t)+1$ has a positive index and decreases by one with
rate $\bar p$ if the particle at $X_0(t)-1$ has a positive index. Thus we
have
%
%
\begin{equation}\label{eq:dR}
\frac{d}{dt}R(t)=p r_0^+(t)-\bar p r_0^-(t).
\end{equation}

Consider the set $A$ of $i$ with a higher particle to $i$'s right, and
the set $B = \{n \dvtx Y_n(t)<Y_{n+1}(t)\}$. By translation invariance, the
density of $A$ is $r^+_0(t)$, and the density of $B$ is $Q(t)$. There is
a bijection between the sets, mapping $i\in A$ to $X_i(t)\in B$. Applying
the mass transport principle (see, e.g.,~\cite{LPbook}), to the
transportation of a unit mass from each $i\in A$ to $X_i(t)\in B$ we find
that $r^+_0(t)=Q(t)$. The same argument shows $r^-_0(t)=\bar{Q(t)}$.
\end{pf}
\begin{pf*}{Proof of Theorem~\ref{T:ASEP-swap}}
Combining the previous three lemmas gives that
\[
\rho/3 = \lim_{t\to\infty} p + Q(t) - 1
\]
(where the limit exists due to the monotonicity proved in
Lemma~\ref{L:Q_mon}). Hence $\lim_{t\to\infty} Q(t) = \frac{2-p}{3}$.
\end{pf*}

\subsection{Joint density}

Throughout this subsection we assume Conjecture~\ref{C:ASEP_lln}.
Under this assumption we can talk about the eventual speed of a particle,
and we know that for large $t$ the empiric speed approximates the eventual
speed. We consider the quantity
\begin{eqnarray*}
R_{a,b}(t) &=& \E\Biggl[ \sum_{j=1}^{\infty} \ind\{ U_j<a,
X_0(t)>X_j(t)\} \cdot\ind[U_0>b] \Biggr],
\\
R_{a,b}(t) &=& \sum_{j>0} \P\bigl(U_0 > b, U_j<a, X_0(t)>X_j(t)\bigr).
\end{eqnarray*}
Thus we ask for $0$ to have speed at least $b$ and count particles of speed
at most $a$ that it overtakes by time $t$. This is of interest for any pair
$-\rho<a<b<\rho$.

\begin{lemma}\label{lemma9.6}
Assume Conjecture~\ref{C:ASEP_lln} holds. Then
\[
R_{a,b}(t) \sim t \int_{-\rho}^a \int_{b}^\rho\frac{y-x}{4\rho^2} \,dy
\,dx = \frac{(\rho+a)(\rho-b)(2\rho+b-a)}{8\rho^2} t.
\]
\end{lemma}

\textit{Note}: this essentially says that the contribution to $R_{a,b}$ from $0$
having speed $y$ (or in $dy$) and $j$'s that have speed $x$ is roughly
$\frac{y-x}{4\rho^2} t$.
\begin{pf*}{Proof of Lemma~\ref{lemma9.6}}
Each particle moves at rate at most $1$, so we have $\P(X_0(t)>X_j(t)) <
\P(\operatorname{Poi}(2t)\geq j)$. This implies that
\[
R_{a,b}(t) = o(1) + \sum_{j=1}^{3t} \P\bigl(U_0 > b, U_j<a,
X_0(t)>X_j(t)\bigr).
\]
The probability that any particle deviates at time $t$ by more than
$\eps$ from its eventual speed is $o(1)$. It follows that
\[
R_{a,b}(t) = o(t) + \sum_{j=1}^{3t} \P\bigl(U_0(t) > b, U_j(t)<a,
X_0(t)>X_j(t)\bigr).
\]
From here on we argue as in the proof of Lemma~\ref{L:R_t_lim}. The
hydrodynamic limit shows that $R_{a,b}(t)$ is asymptotically close to
what it would be if the speeds were independent uniform on
$[-\rho,\rho]$
\begin{eqnarray*}
\frac1t R_{a,b}(t) &=& o(1) + \E\bigl[ \ind\{U_0(t)>b\} \cdot
\nu_t \bigl\{ (x,y)\dvtx x\in\bigl(0,U_0(t)-y\bigr), y<a, x<3 \bigr\} \bigr] \\
&=& o(1) + \frac1{2\rho} \E\bigl[ \ind\{U_0(t)>b\} \\
&&\hspace*{54.7pt}{}\times
\Leb\bigl\{ (x,y)\dvtx x\in\bigl(0,U_0(t)-y\bigr), -\rho\le y<a \bigr\} \bigr]\\
&=& o(1) + \frac1{4\rho^2} \E[ \ind\{U_0(t)>b\}\cdot(2U_0+\rho-a)
(a+\rho)].
\end{eqnarray*}
Simple integration completes the proof.
\end{pf*}

Let $Q_{a,b}(t)$ be the probability of having at time $t$, in positions $0,1$
two particles of speeds in $[b,1]$ and $[-1,a]$, respectively,
\[
Q_{a,b}(t) = \P\bigl(U_{Y_{0}(t)}>b \mbox{ and } U_{Y_{1}(t)}<a\bigr).
\]
We also let $\widetilde Q_{a,b}(t)$ be the probability of having the same
speeds but exchanged
\[
\widetilde Q_{a,b}(t) = \P\bigl(U_{Y_{0}(t)} <a \mbox{ and } U_{Y_{1}(t)} >b\bigr).
\]

\begin{lemma}
Assume Conjecture~\ref{C:ASEP_lln} holds. Then for any $a,b,t$
\[
\frac{d}{dt} R_{a,b}(t) = \bigl(p Q_{a,b}(t) - \bar{p} \widetilde Q_{a,b}(t)\bigr).
\]
\end{lemma}
\begin{pf}
This is an analogue of Lemma~\ref{L:R_t_deriv}. $R_{a,b}(t)$ is the expected
size of the set of $j$'s that are swapped with $0$ at time $t$ (with some
constraints on $U_0,U_j$). This set increases when $0$ has speed at least
$b$ and swaps with a particle of speed at most $a$. Using ergodicity and
translation invariance, just as in Lemma~\ref{L:R_t_deriv}, we find
that the
expected rate at which $j$'s are added to the set is $p Q_{a,b}(t)$.
Similarly, the expected rate at which elements are removed from the set
is $\bar{p} \widetilde Q_{a,b}(t)$. The claim follows.
\end{pf}

Recall that we denote by $\muI$ the joint distribution of $U_0,U_1$ which
we assume exists.
\begin{lemma}
Assume Conjecture~\ref{C:ASEP_lln} holds. Then
\begin{eqnarray*}
\lim_{t\to\infty} Q_{a,b}(t) &=& \muI(U_0<-b \mbox{ and } U_1 >-a) , \\
\lim_{t\to\infty} \widetilde Q_{a,b}(t)
&=& \muI(U_1<-b \mbox{ and } U_0 > -a).
\end{eqnarray*}
\end{lemma}
\begin{pf}
Using $A\approx B$ for $A-B\xrightarrow{t\to\infty}0$, we have
\begin{eqnarray*}
Q_{a,b}(t)
&=& \P\bigl(U_{Y_0(t)}>b \mbox{ and } U_{Y_1(t)}<a\bigr) \\
&\approx& \P\bigl(U_{Y_0(t)}(t) > b \mbox{ and } U_{Y_1(t)}(t)<a\bigr)
\qquad\mbox{by convergence} \\
&=& \P\bigl(Y_0(t) < -b t \mbox{ and } Y_1(t) > 1 - a t\bigr)
\qquad\mbox{since $X_{Y_j(t)}(t)=j$}\\
&=& \P\bigl(X_0(t) < -b t \mbox{ and } X_1(t) > 1 - a t\bigr)
\qquad\mbox{by symmetry} \\
&=& \P\bigl(U_0(t) < -b \mbox{ and } U_1(t) > -a\bigr) \qquad\mbox{by definition} \\
&\approx& \P(U_0 < -b \mbox{ and } U_1 > -a) \qquad\mbox{by convergence,}
\end{eqnarray*}
$\widetilde Q$ is dealt with similarly.
\end{pf}
\begin{pf*}{Proof of Theorem~\ref{T:ASEP-joint}}
Combining the above lemmas and taking the limit as $t\to\infty$ we find
that
\begin{eqnarray*}
\int_{-\rho}^a \int_{b}^\rho\frac{y-x}{4\rho^2} \,d y \,d x
&=& p \muI(U_0<-b, U_1 >-a) - \bar{p}\muI(U_1<-b, U_0 >
-a) \\
&=& \bigl(p\muI- \bar{p}\muII\bigr)(A),
\end{eqnarray*}
where $A=[-\rho,-b)\times(-a,\rho]$. These rectangles determine the
measure $p \muI- \bar{p}\muII$ in the set $\{(x,y)\dvtx-\rho\le x <
y\le
\rho\}$, and differentiating with respect to $a$ and $b$ gives the
statement of the theorem.
\end{pf*}

\subsection{Equal speeds imply interaction}

\mbox{}

\begin{pf*}{Proof of Theorem~\ref{T:together_forever}}
Since we have that $\{J_{0,1}=\infty\} \subset\{U_0=U_1\}$, it suffices
to to prove that $\P(U_0=U_1, J_{0,1}<\infty) = 0$.

In the case of the TASEP the proof is very simple. From
Theorem~\ref{T:ASEP-swap} we know that the probability that particles
$0$ and
$1$ never swap is $1/3$. On the other hand, Theorem~\ref{T:joint2}
implies that $\P(U_0<U_1)=1/3$, and clearly on this event they never
swap. Thus $\P(\mathrm{swap}|U_0\geq U_1)=1$, and the result follows.

The argument for the ASEP mirrors the above, but is more delicate.
Theorem~\ref{T:ASEP-joint} takes on the role of Theorem \ref
{T:joint2}. Start with
%
%
\begin{eqnarray}\label{eq:together1}\quad
\frac{2-p}{3} &=& \lim Q(t)
= \lim_{t\to\infty} \P\bigl(X_0(t)<X_1(t)\bigr) \nonumber\\
&=& \lim_{t\to\infty} \P\bigl(X_0(t)<X_1(t), J_{0,1}<\infty\bigr) +
\P\bigl(X_0(t)<X_1(t), J=\infty\bigr)\\
&=& \P\bigl(\mbox{eventually }X_0(t)<X_1(t)\bigr) + \bar p \P(J_{0,1}=\infty)
.\nonumber
\end{eqnarray}
We also have
%
%
\begin{eqnarray}\label{eq:together2}
&&\P\bigl(\mbox{eventually }X_0(t)<X_1(t)\bigr)\nonumber\\[-8pt]\\[-8pt]
&&\qquad= \P(U_0<U_1)
+ \E\bigl[ \ind[U_0=U_1] \ind[J_{0,1}<\infty] ( \bar p + p
e^{-J_{0,1}} ) \bigr].\nonumber
\end{eqnarray}
[Compare with \eqref{eq:J} and the discussion around it.] Combining
\eqref{eq:together1} and \eqref{eq:together2} and noting that
$\P(J_{0,1}=\infty)=\P(J_{0,1}=\infty, U_0=U_1)$ we get
%
%
\begin{eqnarray}\label{eq:together3}
\frac{2-p}{3} &=& \P(U_0<U_1) + \bar p \P(U_0=U_1) \nonumber\\[-8pt]\\[-8pt]
&&{}+
\E\bigl[ \ind[U_0=U_1] \ind[J_{0,1}<\infty] p
e^{-J_{0,1}}\bigr].\nonumber
\end{eqnarray}
On the other hand, integrating Theorem~\ref{T:ASEP-joint} gives
\[
\frac{2p-1}{3} = p \P(U_0<U_1) - \bar p \P(U_0>U_1),
\]
which implies
\[
\frac{2-p}{3} = \P(U_0<U_1) + \bar p \P(U_0=U_1).
\]
Together with \eqref{eq:together3} this implies
\[
\E\bigl[ \ind[U_0=U_1] \ind[J_{0,1}<\infty] p e^{-J_{0,1}}\bigr] = 0,
\]
and so $\P(U_0=U_1, J_{0,1}<\infty)=0$ as needed.
\end{pf*}

This can be extended to other particles with equal speeds. Let
$J_{i,j}$ be
the total time that particles $i$ and $j$ are in adjacent positions.
\begin{lemma}\label{L:diagswap}
For any $k>i$, a.s.
\[
k=\min\{j>i \dvtx U_j=U_i\} \quad\implies\quad J_{i,k}=\infty.
\]
Consequently, in the TASEP every two particles in the same convoy swap
eventually.
\end{lemma}
\begin{pf}
Clearly this only depends on $k-i$. We proceed by induction on $k-i$. For
$k=i+1$ this is just Theorem~\ref{T:together_forever}. The key to the induction
step is to show that if $U_0\neq U_1$ then there is a transformation of
the probability space that swaps the eventual trajectories of 0 and 1
(and hence their speeds), keeps all other trajectories the same, and has
finite Radon--Nikodym derivative. It follows that applying this
transformation results in an absolutely continuous measure for the
trajectories. If we assume the lemma for $k$ and $1$, then
\[
\P(k=\min\{j>1 \dvtx U_j=U_1\} \mbox{ and } J_{1,k}<\infty) = 0,
\]
and hence by absolute continuity the result holds for $k,0$.

Recall the $\sigma$-field $\F_{0,1}$ of the trajectories of all particles
except $0$ and $1$. If $U_0>U_1$ the transformation just eliminates all
interactions between $0$ and $1$. This has the effect of exchanging their
trajectories from some point on. Given $\F_{0,1}$, the probability of no
interaction between $0$ and $1$ is $e^{-J_{0,1}}$. The Radon--Nikodym
derivative is at most $e^{J_{0,1}}<\infty$ (on $U_0\neq U_1$).

If $U_0<U_1$ we define the transformation as follows: consider the first
time $\tau$ at which either $0$ or $1$ swaps with some other particle,
and replace all interactions between $0$ and $1$ by a unique interaction
between $0$ and $1$ at a time uniform on $[0,\tau]$. In the ASEP, we make
this new interaction exchange $0$ and $1$. The probability of this
pattern of interactions between $0$ and $1$, given $\F_{0,1}$ is $p\tau
e^{-J_{0,1}}$, thus the Radon--Nikodym derivative in this case is at most
$e^{J_{0,1}}/(p\tau)<\infty$.

Finally, in the TASEP, since any pair of consecutive particles in a
convoy a.s. swap and particles never unswap, it follows that all pairs
eventually swap.
\end{pf}

\section*{Acknowledgments}
The authors wish to thank James Martin, Pablo Ferrari and B\'{a}lint
Vir\'{a}g
for useful discussions.

%
%

%
\printaddresses

\end{document}